\documentclass[12pt]{article}
\usepackage{makeidx}
\usepackage{amsfonts}
\usepackage{amssymb}
\usepackage{amsmath,amsthm}
\usepackage{graphicx}
\usepackage{color}
\usepackage{ulem}
\makeatletter
\@addtoreset{equation}{section}
\@addtoreset{section}{part}
\makeatother

\newtheorem{theorem}{Theorem}[section]
\newtheorem{lemma}[theorem]{Lemma}
\newtheorem{corollary}[theorem]{Corollary}
\newtheorem{proposition}[theorem]{Proposition}

\newtheorem{remark}[theorem]{Remark}

\newtheorem{definition}[theorem]{Definition}

\renewcommand{\P} {\mathbb{P}}
\newcommand{\R} {\mathbb{R}}

\begin{document}

\title{Uniqueness of the solution of the filtering equations in spaces of
measures}
\author{Dan Crisan, Etienne Pardoux}
\maketitle

\begin{abstract}
Nonlinear filtering is a pivotal problem that has attracted significant attention from mathematicians, statisticians, engineers, and various other scientific disciplines. The solution to this problem is governed by the so-called "filtering equations." In this paper, we investigate the uniqueness of solutions to these equations within measure spaces and introduce a novel, generalized framework for this analysis. Our approach provides new insights and extends the applicability of existing theories in the study of nonlinear filtering.

\end{abstract}

\section{Introduction}

The objective of nonlinear filtering is to estimate an evolving dynamical system, modeled by a stochastic process $X$, referred to as the signal process. This signal process cannot be directly observed; instead, it is inferred through a related process $Y$, known as the observation process. The filtering problem involves determining the conditional distribution of the signal at the current time, based on the observation data accumulated up to that time. For a comprehensive treatment of the filtering problem and a historical overview of its extensive development over the past 80 years, starting with the foundational work of Kolmogorov, Krein, and Wiener, readers are encouraged to consult \cite{cr}.

The conditional distribution of the signal  given the observation is the solution of a nonlinear
evolution equation. Moreover, it has a version which satisfies a linear
evolution equations. In what follows we shall refer to these to equations as
the filtering equations. The filtering equations have been studied at length
by many different contributors to the filtering problem and in different
frameworks. The most popular continuous time framework is that where the
signal $X$ satisfies a stochastic differential equation driven by a
multi-dimensional Brownian motion denoted by $V$, and the observation $Y$
satisfies an evolution equation of the form 
\begin{equation*}
dY_{t}=h(t,X_{t})\,\mathrm{d}t+d{W}_{t},
\end{equation*}%
where the driving Brownian motion $W$ is independent of $X$. In this case,
the filtering equations are well understood and their solution is shown to
be unique in suitably chosen spaces of measures, see, e.g., \cite{bc}, for
details.

In this set-up, the equations satisfied by the signal and the
observation are asymmetric. There is no dependence of $X$ on $Y$, and the
observational noise $W$ is chosen independent of the signal and with a
coefficient that does not depend on either $X$ or $Y$. Possibly a symmetric
set-up would require the pair ${\mathcal{X}}:=(X,Y)$ to satisfy a stochastic
differential of the form 
\begin{equation*}
d{\mathcal{X}}_{t}=a(t,{\mathcal{X}}_{t})\,\mathrm{d}t+b(t,{\mathcal{X}}%
_{t})d{B}_{t},
\end{equation*}%
where $B:=(V,W)$. However, perfect symmetry is not  really possible. For example, should
the diffusion observation coefficient depend on the signal, the solution
will degenerate, see, e.g., \cite{ckx} for details.

In this paper, we
introduce a framework which is as close as possible to a symmetric one, see
equations (\ref{X}) and (\ref{Y}) below. This new framework represents a
generalization of several existing ones, including those covered in \cite%
{bensoussan} and \cite{pp2}. A simplified version of the present set-up
(where the signal and the observation noises are independent) was studied in 
\cite{FDP}.

In the new framework the coefficients of the stochastic differential
equation satisfied by $X$ depend on the pair $(X,Y)$. Moreover $X$ is driven
by the pair $(V,W)$ and not just by $V$. Put it differently the Brownian
motion driving the signal is corellated with the Brownian motion driving the
observation. Moreover we do not assume that the diffusion coefficient in the
observation equation is invertible. In the degenerate case when the
diffusion coefficient is equal to 0, the observation becomes independent of
the signal: in this case the conditional distribution of the signal
coincides with the (prior) distribution of the signal.

Within this framework we obtain the filtering equations and show that they
have a unique solution within a suitably chosen space of measures. The following are  contributions of this paper:

\begin{itemize}
\item We deduce the filtering equations for the signal and the observation
processes satisfying equations (\ref{X})--(\ref{Y}) below. The coefficients
of both equations depend on the signal-observation pair. Moreover, as
opposed to frameworks treated elsewhere we do not assume an invertible
diffusion term in the observation equation. We do this at the expense of
assuming a special form for the observation function $h$, see (\ref{h})
below.

\item We show the equivalence of the uniqueness properties of
the two filtering equations, namely the nonlinear Kushner--Stratonovich
equation for the conditional distribution, and the linear Zakai equation for
the ``unnormalized conditional distribution'', see equations %
\eqref{eq:filterEqns:ks} and, respectively, %
\eqref{eq:filterEqns:zakaifunctint} below.

\item We establish the uniqueness of the solution of the equation satisfied
by the unnormalized conditional distribution of the signal in the space of
measure valued processes.
\end{itemize}

 Let us now explain the novelty of our proof of the uniqueness
of the Zakai equation. The so-called ``duality argument'' is a standard
method for proving uniqueness of the solution of a linear equation.
Bensoussan applied the duality argument in \cite{bensoussan} to show the
uniqueness of the solution of the Zakai equation.  In the framework treated
in \cite{bensoussan}, the measure valued solution of the Zakai equation,
which is a stochastic partial differential equation was paired with the
(function valued) solution of a deterministic backward PDE.\footnote{%
The duality property is shown through the use of a collection of
exponential martingales first introduced in the filtering framework by
Krylov and Rozovskii.\cite{KR}} The proof in \cite{bensoussan} involves the
use of an It\^o's formula to deduce the evolution of the solution of a
deterministic backward PDE integrated with respect to the measure valued process
that is the solution of the Zakai equation. The same argument cannot be
applied in the framework treated here since the coefficients of the
operators appearing in the Zakai equation (hence also in the adjoint
backward PDE) depend upon the current observation. As a result, the dual of
the solution of the Zakai equation would now be the solution of a backward
PDE starting at the time $t$ and run backwards in the interval $[0,t]$ and,
at any time $s\in (0,t)$, would be a function of the observations on the
interval $[s,t]$, making it random and, moreover, anticipating at any time $%
s\in [0,t]$ the future of the observation process. As a result, the It\^o
formula can no longer be applied. In order to circumvent this difficulty, we replace
the dual deterministic PDE by a \emph{backward} Stochastic Partial Differential Equation''
(or BSPDE for short). For this, we exploit recent results from Du, Meng \cite%
{DM} and Du, Tang , Zhang \cite{DTZ} which we need to adapt to our framework
which requires a complex--valued BSPDE, equivalent to a system of two real
valued BSPDEs, and construct Sobolev space valued solutions of our system of
BSPDEs. Moreover, we exploit classical Sobolev embedding theorems, in order
to deduce enough smoothness of the solution of the BSPDE so that it
can be integrated against an arbitrary measure--valued solution to the
Zakai equation. 

\bigskip

The paper is structured as follows. In Section 2 we introduce the filtering
framework and the two sets of assumptions that ensure the existence and
uniqueness of the filtering equations (see \textbf{Assumption E} and \textbf{%
Asssumption U} below). In Section 3 we deduce the Kallianpur-Striebel
formula (Proposition \ref{prop:filterEqns:kallstrie}) which implies the
existence of an unnormalized version of the conditional distribution of the
signal. Next we deduce the filtering equation (\ref{eq:filterEqns:ks}) for
the conditional distribution of the signal,  which is a non linear SPDE, as well as  a linear equation
(\ref{eq:filterEqns:zakaifunctint}) for the unnormalized version, see Theorem \ref%
{th:ex} below. Finally we show in Theorem \ref{equivuniq} that equation (\ref%
{eq:filterEqns:ks}) has a unique solution if and only if (\ref%
{eq:filterEqns:zakaifunctint}) has a unique solution. In Section 4 we
introduce several results pertaining to a class of backward stochastic
partial differential equations that are used in the subsequent section as
well as a useful It\^o type formula, see Theorem \ref{Zakaiextended} and
some preliminary Lemmas. The paper is concluded with Section 5 where the
uniqueness of the solution of equation (\ref{eq:filterEqns:zakaifunctint}),
respectively (\ref{eq:filterEqns:ks}) is proved, see Theorem \ref{th:Unique}
below. Finally in the Appendix we recall the definition of the Moore--Penrose pseudo--inverse $A^+$
of a (possibly) rectangular matrix $A$, and prove that the mapping  $A\mapsto A^+$ is measurable.

\section{Framework}

Let $(\Omega ,\mathcal{F},\mathbb{P})$ be a probability space together with
a filtration $(\mathcal{F}_{t})_{t\geq 0}$ which satisfies the usual
conditions\footnote{%
The probability space $(\Omega ,\mathcal{F},\mathbb{P})$ together with the
filtration $(\mathcal{F}_{t})_{t\geq 0}$ satisfies the usual conditions
provided: a. $\mathcal{F}$ is complete i.e. $A\subset B$, $B\in \mathcal{F}$
and $\mathbb{P}(B)=0$ implies that $A\in \mathcal{F}$ and $\mathbb{P}(A)=0$,
b. The filtration $\mathcal{F}_{t}$ is right continuous i.e. $\mathcal{F}
_{t}=\mathcal{F}_{t+}$, and,  c. $\mathcal{F}_{0}$ (and consequently all $\mathcal{%
\ F}_{t}$ for $t\geq 0$) contains all the $\mathbb{P}$-null sets.}. 
We recall that to such a filtration we associate the $\sigma$--algebra $\mathcal{P}(\mathcal{F}_t)$ of progressively measurable subsets of 
$\Omega\times\R_+$, which is the class of sets $A\subset\Omega\times\R_+$ which are such that for all $t\ge0$,
\[ A\cap(\Omega\times[0,t])\in\mathcal{F}_t\otimes\mathcal{B}_{[0,t]},\]
where $\mathcal{B}_{[0,t]}$ denotes the $\sigma$--algebra of Borel measurable subsets of $[0,t]$.
On $(\Omega ,\mathcal{F},\mathbb{P})$ we consider a $\mathcal{P}(\mathcal{F}_{t})$-measurable\footnote{This means that the mapping 
$(\omega,t)\mapsto(X_t(\omega),Y_t(\omega))$ is $\mathcal{P}(\mathcal{F}_{t})$-measurable.}
process $(X,Y)$ with continuous paths. The process $X$ is called the \emph{\
signal} process and is assumed to take values in $\mathbb{R}^{d}$ (termed as the state
space). The process $Y$ is assumed to take values in $\mathbb{R}^{d^{\prime
}}$ and is called the \emph{observation} process.

We will assume that the processes $(X,Y)$ satisfy the following systems of
stochastic differential equations 
\begin{eqnarray}
X_{t} &=&X_{0}+\int_{0}^{t}f(s,X_{s},Y_{s})\,\mathrm{d}s+%
\int_{0}^{t}g(s,X_{s},Y_{s})\,\mathrm{d}V_{s}+\int_{0}^{t}\bar{g}%
(s,X_{s},Y_{s})\,\mathrm{d}W_{s},  \label{X} \\
Y_{t} &=&Y_{0}+\int_{0}^{t}h(s,X_{s},Y_{s})\,\mathrm{d}s+%
\int_{0}^{t}k(s,Y_{s})d{W}_{s},  \label{Y}
\end{eqnarray}%
where $V$ and $W$ are mutually independent $\ell $
(resp. $\ell ^{\prime }$) dimensional standard Brownian motions, and $f$,$g,\bar{g},h,k$ satisfy
suitable conditions so that the system (\ref{X})+(\ref{Y}) has a unique
solution (see \textbf{Assumption E} below). In addition, we assume that 
\begin{equation}
h(s,x,y)=h_{1}(s,y)+k(s,y)h_{2}(s,x,y).  \label{h}
\end{equation}%
Let $\mathcal{B}(\mathbb{R}^{d})$ and $\mathcal{B}(\mathbb{R}^{d}\times \mathbb{R}^{d^{\prime }})$ be the
associated product Borel $\sigma $-algebra on $\mathbb{R}^{d}$ and, respectively, 
$\mathbb{R}^{d}\times \mathbb{R%
}^{d^{\prime }}$  and let  $b\mathcal{B}(\mathbb{R}^{d})$ and 
$b\mathcal{B}(\mathbb{R}^{d}\times \mathbb{R}%
^{d^{\prime }})$ be the space of bounded 
$\mathcal{B}(\mathbb{R}^{d}\times 
\mathbb{R}^{d^{\prime }})$, respectively, 
$\mathcal{B}(\mathbb{R}^{d}\times 
\mathbb{R}^{d^{\prime }})$ measurable functions. Let $A_{s}$ be the
following differential operator%
\begin{eqnarray*}
A_{s}\varphi \left( x\right)  &=&\sum_{i=1}^{d}f^{i}\left( s,x,Y_{s}\right)
\partial _{i}\varphi \left( x\right) +\frac{1}{2}\sum_{i,j=1}^{d}a^{ij}%
\left( s,x,Y_{s}\right) \partial _{i}\partial _{j}\varphi \left( x\right),  \\
a^{ij}\left( s,x,Y_{s}\right)  &=&\sum_{p=1}^{\ell }g^{ip}g^{jp}\left(
s,x,Y_{s}\right) +\sum_{p=1}^{\ell ^{\prime }}\bar{g}^{ip}\bar{g}^{jp}\left(
s,x,Y_{s}\right) .
\end{eqnarray*}%
We will impose the following sets of assumptions on the coefficients of the
system (\ref{X})+(\ref{Y}):

\textbf{Assumption E}. The functions%
\begin{eqnarray*}
f &:&[0,\infty )\times \mathbb{R}^{d}\times \mathbb{R}^{d^{\prime
}}\rightarrow \mathbb{R}^{d} \\
g &:&[0,\infty )\times \mathbb{R}^{d}\times \mathbb{R}^{d^{\prime
}}\rightarrow \mathbb{R}^{d\times \ell } \\
\bar{g} &:&[0,\infty )\times \mathbb{R}^{d}\times \mathbb{R}^{d^{\prime
}}\rightarrow \mathbb{R}^{d\times \ell ^{\prime }} \\
h &:&[0,\infty )\times \mathbb{R}^{d}\times \mathbb{R}^{d^{\prime
}}\rightarrow \mathbb{R}^{d^{\prime }} \\
h_{2} &:&[0,\infty )\times \mathbb{R}^{d}\times \mathbb{R}^{d^{\prime
}}\rightarrow \mathbb{R}^{\ell^{\prime }} \\
h_{1} &:&[0,\infty )\times \mathbb{R}^{d^{\prime }}\rightarrow \mathbb{R}%
^{d^{\prime }} \\
k &:&[0,\infty )\times \mathbb{R}^{d^{\prime }}\rightarrow \mathbb{R}%
^{d^{\prime }\times \ell^{\prime }}
\end{eqnarray*}%
have the following properties:

\begin{itemize}
\item $f,g,\bar{g},h,h_{1}$ and $k$ are locally Lipschitz in the $(x,y)$ variables. In other words, for any $R>0$, we have that  
\begin{equation*}
\left\vert \left\vert f\left( t,x_{1},y_{1}\right) -f\left(
t,x_{2},y_{2}\right) \right\vert \right\vert \leq K_R\left( \left\vert
\left\vert x_{1}-x_{2}\right\vert \right\vert +\left\vert \left\vert
y_{1}-y_{2}\right\vert \right\vert \right) ,~~x_{1},x_{2}\in \mathbf{B}%
_{R}^{d},~~y_{1},y_{2}\in \mathbf{B}_{R}^{d^{\prime }},
\end{equation*}%
where $\mathbf{B}_{R}^{d}$ (resp. $\mathbf{B}_{R}^{d^{\prime }}$) is the
ball of centre $0$ and radius $R$ in $\mathbb{R}^{d}$ (resp. in $\mathbb{R}^{d^{\prime }}$) and $K_R$ is a constant  which may depend upon $R$, but is    independent of all
variables, with a similar condition imposed on $g,\bar{g},h,h_{1}$ and $k.$

\item $f,g,\bar{g},h,h_{1},h_{2}$ and $k$ satisfy a linear growth condition  in the $(x,y)$ variables.
In other words,%
\begin{equation*}
\left\vert \left\vert f\left( t,x,y\right) \right\vert \right\vert \leq
K\left( 1+\left\vert \left\vert x\right\vert \right\vert +\left\vert
\left\vert y\right\vert \right\vert \right) ~~x_{1},x_{2}\in \mathbb{R}%
^{d},~~y_{1},y_{2}\in \mathbb{R}^{d^{\prime }},
\end{equation*}%
where $K$ is a constant independent of all variables, with a similar
condition imposed on $g,\bar{g},h,h_{1}$ and $k.$
\end{itemize}

We also assume that $X_{0}$ and $Y_{0}$ have finite second moments, that is $\mathbb{%
E}\left[ || X_{0} || ^{2}+|| Y_{0} || ^{2}\right] <\infty $. The following is a classical result in the theory of stochastic differential equations, see e.g., \cite{ks} for a proof.   

\begin{remark}
Under \textbf{Assumption E}, the system (\ref{X})+(\ref{Y}) has a unique
global solution. Moreover, for any $T>0$, we have that 
\begin{equation}
\mathbb{E}\left[ \sup_{s\in \left[ 0,T\right] }\left\vert \left\vert
X_{s}\right\vert \right\vert ^{2}\right] +\mathbb{E}\left[ \sup_{s\in \left[
0,T\right] }\left\vert \left\vert Y_{s}\right\vert \right\vert ^{2}\right]
<\infty .  \label{squareint}
\end{equation}
\end{remark}

\textbf{Assumption U}. The functions $f$, $g,\bar{g},h_2$ are bounded on $%
[0,T]\times \mathbb{R}^{d}\times \mathbb{R}^{d^{\prime }}$ for arbitrary\ $%
T>0$. The functions $h_{1}$ and $k$ are bounded on $[0,T]\times \mathbb{R}%
^{d^{\prime }}$ for arbitrary\ $T>0$. Moreover, for some integer $n>%
\frac{d}{2}+2$, all the partial derivatives of the functions $f$, $%
g,\bar{g},h$ in the $x$ variable with multi-index $\alpha $, such that $%
|\alpha |\leq n$, are bounded on $[0,T]\times \mathbb{R}^{d}\times \mathbb{R}%
^{d^{\prime }}$ for arbitrary\ $T>0$ 

\begin{remark}
In the following, the evolution equation of the conditional distribution of
the signal  given the observation  
 will be derived under \textbf{Assumption E}, whilst the
uniqueness of its solution will be proved under the joint \textbf{%
Assumptions E+U}.
\end{remark}

Let $\{\mathcal{Y}_{t},\ t\geq 0\}$ be the usual augmentation of the
filtration associated with the process $Y$, viz 
\begin{equation}
\mathcal{Y}_{t}=\bigcap_{\varepsilon >0}\sigma (Y_{s},\ s\in \lbrack
0,t+\varepsilon ])\vee \mathcal{N},\ \ \ \mathcal{Y}=\bigvee_{t\in \mathbb{R}%
_{+}}\mathcal{Y}_{t}.  \label{yt}
\end{equation}%
where $\mathcal{N}$ is the class of all $\mathbb{P}$-null sets. Note that $%
Y$ is $\mathcal{P}(\mathcal{F}_{t})$-measurable, hence $\mathcal{Y}_{t}\subset \mathcal{F}_{t}$ for all $t\ge0$.\ 

\begin{definition}
The filtering problem consists in determining the conditional distribution $%
\varsigma _{t}$ of the signal $X$ at time $t$ given the information
accumulated from observing $Y$ in the interval $[0,t]$; that is, for any $%
\varphi \in b\mathcal{B}(\mathbb{R}^d)$, 
\begin{equation}
\varsigma _{t}(\varphi )=\mathbb{E}[\varphi (X_{t})\mid \mathcal{Y}_{t}].
\label{nfp}
\end{equation}
\end{definition}

\section{The Filtering Equations}

In the following we deduce the evolution equation for $\varsigma_t $. Let $%
\tilde{W}=\{\tilde{W}_{t},\ t\geq 0\}$ be the process defined as%
\begin{equation*}
\tilde{W}_{t}=W_{t}+\int_{0}^{t}h_{2}(s,X_s,Y_s)ds.
\end{equation*}%
We shall construct a new measure under which $\tilde{W}$ becomes a Brownian
motion and $\varsigma $ has a representation in terms of an associated
unnormalised version $\pi $. This $\pi $ is then shown to satisfy a linear
evolution equation which leads to the evolution equation for $\varsigma $ by
an application of It\^{o}'s formula. Define $Z=\left( Z_{t}\right) _{t\geq
0} $ to be the exponential local martingale 
\begin{equation}
Z_{t}=\exp \left( -\int_{0}^{t}h_{2}\left( s,X_{s},Y_{s}\right) ^{\top
}dW_{s}-\frac{1}{2}\int_{0}^{t}\left\vert h_{2}\left( s,X_{s},Y_{s}\right)
\right\vert ^{2}ds\right) .  \label{zt}
\end{equation}%
We will work under the following additional assumption:

\textbf{Assumption M.} We assume that 
\begin{equation}
\mathbb{E}\left[ Z_{t}\right] =1,\quad \forall t>0.  \label{condonh}
\end{equation}%
In particular this implies that $Z$ is a genuine martingale (not just a
local martingale).

\begin{remark}
 There are several assumptions under which (\ref{condonh})
holds, see e.g. Proposition 2.50 in \cite{pardouxrascanu}. The sufficient
condition (ii) of that Proposition requires that for each $t>0$, there
exists $\gamma>0$ such that 
\begin{equation*}
\sup_{0\le s\le t}\mathbb{E}[\exp(\gamma |h_2(s,X_s,Y_s)|^2)]<\infty.
\end{equation*}
This condition is, in particular, satisfied if $h_2(t,X_t,Y_t)$ is a Gaussian
process with locally bounded mean and variance, and also clearly if $h_2$ is bounded. It also follows that Assumption {\bf U} implies 
Assumption {\bf M}. 
\end{remark}

Let $\mathbb{\tilde{P}}$ be the probability measure defined on the field $%
\bigcup_{0\leq t<\infty }\mathcal{F}_{t}$ that is specified by its
Radon--Nikodym derivative $Z_{t}$ on each $\mathcal{F}_{t}$ with respect to
the corresponding trace of $\mathbb{P}$; that is, for each $t\geq 0$:%
\begin{equation*}
\left. \frac{\mathrm{d}\tilde{\mathbb{P}}}{\mathrm{d}\mathbb{P}}\right\vert
_{\mathcal{F}_{t}}=Z_{t}.
\end{equation*}%
$\mathbb{\tilde{P}}$ restricted to each $\mathcal{F}_{t}$ is equivalent to $%
\mathbb{P}$ since $Z_{t}$ is a positive random variable\footnote{%
Note that we have not defined $\mathbb{\tilde{P}}$ on $\mathcal{F}_{\infty }$%
, where $\mathcal{F}_{\infty }=\bigvee_{t=0}^{\infty }\mathcal{F}_{t}=\sigma
\left( \bigcup_{0\leq t<\infty }\mathcal{F}_{t}\right) .$}.

Let $\tilde{Z}=\{\tilde{Z}_{t},\ t\geq 0\}$ be the process defined as $%
\tilde{Z}_{t}=Z_{t}^{-1}$ for $t\geq 0$. Under $\tilde{\mathbb{P}}$, $\tilde{%
Z}_{t}$ satisfies the following stochastic differential equation, 
\begin{equation}
\mathrm{d}\tilde{Z}_{t}=\tilde{Z}_{t}h_{2}(t,X_t,Y_t)^{\top }\,\mathrm{d}%
\tilde{W}_{t}=\sum_{j=1}^{\ell^{\prime }}\tilde{Z}_{t}h_{2}^{j}(t,X_t,Y_t)\,%
\mathrm{d}\tilde{W}_{t}^{j}  \label{ztilde1}
\end{equation}%
and since $\tilde{Z}_{0}=1$, 
\begin{equation}
\tilde{Z}_{t}=\exp \left( \sum_{j=1}^{\ell^{\prime
}}\int_{0}^{t}h_{2}^{j}(s,X_s,Y_s)\, \mathrm{d}\tilde{W}_{s}^{j}-\frac{1}{2}%
\sum_{j=1}^{\ell^{\prime }}\int_{0}^{t}h_{2}^{j}(s,X_s,Y_s)^{2}\,\mathrm{d}%
s\right) ,  \label{ztilde2}
\end{equation}
then $\tilde{\mathbb{E}}[\tilde{Z}_{t}]=\mathbb{E}[\tilde{Z}_{t}Z_{t}]=1$.
So $\tilde{Z}$ is an $\mathcal{F}_{t}$--martingale under $\tilde{%
\mathbb{P}}$ and%
\begin{equation*}
\left. \frac{\mathrm{d}\mathbb{P}}{\mathrm{d}\tilde{\mathbb{P}}}\right\vert
_{\mathcal{F}_{t}}=\tilde{Z}_{t}\quad \mathrm{\ for\ }t\geq 0.
\end{equation*}%
The probability measures $\mathbb{P}$ and $\mathbb{\tilde{P}}$ are therefore equivalent on each $\sigma$-field $%
\mathcal{F}_{t}$ for any $t\geq 0$. The following proposition is a direct
consequence of Girsanov's theorem:

\begin{proposition}
\label{prop:filterEqns:YisBM} If assumption {\bf M} holds,
then under $\tilde{\mathbb{P}}$ the process $\tilde{W}$ is a Brownian motion
independent of $V$.
\end{proposition}

\begin{remark}
Since $\mathbb{P}$ and $\mathbb{\tilde{P}}$ are absolutely continuous with
respect to each other, they have the same class of null sets $\mathcal{N}$
and therefore the (augmented) observation filtration is the same both under $%
\mathbb{P}$ and under $\mathbb{\tilde{P}}$.
\end{remark}

The following proposition is a consequence of the independent increments
property of the process $\tilde{W}$ under $\tilde{\mathbb{P}}$.

\begin{proposition}
\label{prop:filterEqns:p4} Let $U$ be an integrable $\mathcal{F}_{t}$%
-measurable random variable. Then we have 
\begin{equation}
\tilde{\mathbb{E}}[U\mid \mathcal{Y}_{t}]=\tilde{\mathbb{E}}[U\mid \mathcal{Y%
}].  \label{yty}
\end{equation}
\end{proposition}

\begin{proof}
Under \textbf{Assumption E}, the process $Y$ is the unique strong solution
of the equation 
\begin{equation*}
Y_{t}=Y_{0}+\int_{0}^{t}h_{1}(s,Y_{s})\,\mathrm{d}s+\int_{0}^{t}k(s,Y_{s})d{%
\tilde{W}}_{s}
\end{equation*}%
driven by the Brownian motion ${\tilde{W}}$ (under $\tilde{\mathbb{P}})$\
and we deduce from this that 
\begin{equation*}
\mathcal{Y}\subset \mathcal{Y}_{t}\vee \mathcal{F}^{t,\tilde{W}}
\end{equation*}%
where $\mathcal{F}^{t,\tilde{W}}=\sigma \ \left( \tilde{W}_{t+s}-\tilde{W}%
_{t}|\ \ s>0\right) $. Moreover $\mathcal{F}^{t,\tilde{W}}$ is independent of $%
\mathcal{F}_{t}\supseteq \mathcal{Y}_{t}$ under $\tilde{\mathbb{P}}$. It
follows that since $U$ is $\mathcal{F}_{t}$ measurable, 
\begin{equation*}
\tilde{\mathbb{E}}[U\mid \mathcal{Y}]=\tilde{\mathbb{E}}[\tilde{\mathbb{E}}%
[U\mid \mathcal{Y}_{t}\vee \mathcal{F}^{t,\tilde{W}}]\mid \mathcal{Y}]=%
\tilde{\mathbb{E}}[U\mid \mathcal{Y}_{t}].
\end{equation*}
\end{proof}

In the following, the notation $\tilde{\mathbb{P}}(\mathbb{P})$-a.s. means
that the result holds both $\mathbb{\tilde{P}}$-a.s. and $\mathbb{P}$-a.s. 

\begin{proposition}[\textbf{Kallianpur--Striebel}]
\label{prop:filterEqns:kallstrie} If assumption {\bf M} holds, for every $\varphi \in b\mathcal{B}(\mathbb{R}^d)$, and $t\in
\lbrack 0,\infty )$, 
\begin{equation}
\varsigma _{t}(\varphi )=\frac{\tilde{\mathbb{E}}[\tilde{Z}_{t}\varphi
(X_{t})\mid \mathcal{Y}]}{\tilde{\mathbb{E}}[\tilde{Z}_{t}\mid \mathcal{Y}]}%
\quad \quad \tilde{\mathbb{P}}(\mathbb{P})\text{-}\mathrm{a.s.}
\label{kallstrie}
\end{equation}
\end{proposition}

\begin{proof}
Since $\varsigma _{t}(\varphi )=\mathbb{E}[\varphi (X_{t})|\mathcal{Y}_{t}]$%
, for any $A\in \mathcal{Y}_{t}$, 
\begin{equation*}
\mathbb{E}[1_{A}\varsigma _{t}(\varphi )]=\mathbb{E}[1_{A}\varphi (X_{t})]\,.
\end{equation*}
Consequently 
\begin{equation*}
\tilde{\mathbb{E}}[1_{A}\varsigma _{t}(\varphi )\tilde{Z}_{t}]=\tilde{%
\mathbb{E}}[1_{A}\varphi (X_{t})\tilde{Z}_{t}]
\end{equation*}%
\begin{equation*}
\tilde{\mathbb{E}}[\tilde{\mathbb{E}}[1_{A}\varsigma _{t}(\varphi )\tilde{Z}%
_{t}|\mathcal{Y}_{t}]]]=\tilde{\mathbb{E}}[\tilde{\mathbb{E}}[1_{A}\varphi
(X_{t})\tilde{Z}_{t}|\mathcal{Y}_{t}]]]
\end{equation*}%
\begin{equation*}
\tilde{\mathbb{E}}[1_{A}\varsigma _{t}(\varphi )\tilde{\mathbb{E}}[\tilde{Z}%
_{t}\mid \mathcal{Y}_t]]=\tilde{\mathbb{E}}[1_{A}\tilde{\mathbb{E}}[\tilde{Z}%
_{t}\varphi (X_{t})\mid \mathcal{Y}_{t}]].
\end{equation*}
It follows that $\tilde{\mathbb{E}}\left[1_{A}\left\{\varsigma _{t}(\varphi )%
\tilde{\mathbb{E}}[\tilde{Z}_{t}\mid \mathcal{Y}_t]\mathbb{-}\tilde{\mathbb{E%
}}[\tilde{Z}_{t}\varphi (X_{t})\mid \mathcal{Y}_{t}]\right\}\right]=0$ for
any $A\in \mathcal{Y}_{t}$. Since 
\begin{equation*}
\varsigma _{t}(\varphi )\tilde{\mathbb{E}}[\tilde{Z}_{t} \mathcal{Y}_t] 
\mathbb{-}\tilde{\mathbb{E}}[\tilde{Z}_{t}\varphi (X_{t}) |\mathcal{Y}_t]
\end{equation*}
is $\mathcal{Y}_{t}$-measurable, we deduce that 
\begin{equation*}
\varsigma _{t}(\varphi )\tilde{\mathbb{E}}[\tilde{Z}_{t}\mid \mathcal{Y}_t]%
\mathbb{-}\tilde{\mathbb{E}}[\tilde{Z}_{t}\varphi (X_{t})\mid \mathcal{Y}%
_{t}]=0
\end{equation*}%
$\tilde{\mathbb{P}}$- almost surely. Hence (\ref{kallstrie}) holds by
Proposition \ref{prop:filterEqns:p4}.
\end{proof}

Let $\zeta =\{\zeta _{t},\ t\geq 0\}$ be the process defined by 
\begin{equation}
\zeta _{t}=\tilde{\mathbb{E}}[\tilde{Z}_{t}\mid \mathcal{Y}],
\label{eq:filterEqns:zeta}
\end{equation}%
then as $\tilde{Z}_{t}$ is an $\mathcal{F}_{t}$-martingale under $\tilde{%
\mathbb{P}}$ and $\mathcal{Y}_{s}\subseteq \mathcal{F}_{s}$, it follows that
for $0\leq s<t$, 
\begin{equation*}
\tilde{\mathbb{E}}[\zeta _{t}\mid \mathcal{Y}_{s}]=\tilde{\mathbb{E}}[\tilde{%
Z}_{t}|\mathcal{Y}_{s}] =\tilde{\mathbb{E}}\left[ \tilde{\mathbb{E}}[\tilde{Z%
}_{t}\mid \mathcal{F}_{s}]\mid \mathcal{Y}_{s}\right] =\tilde{\mathbb{E}}[%
\tilde{Z}_{s}\mid \mathcal{Y}_{s}]= \tilde{\mathbb{E}}[\tilde{Z}_{s}\mid 
\mathcal{Y}]=\zeta _{s},
\end{equation*}
where the penultimate equality follows by Proposition \ref%
{prop:filterEqns:p4}. Therefore by Doob's regularization theorem (see e.g.
Theorem 3.13 in \cite{ks}), since the filtration $\mathcal{Y}_{t}$ satisfies
the usual conditions we can choose a c\`{a}dl\`{a}g version of $\zeta _{t}$
which is a $\mathcal{Y}_{t}$-martingale.
In what follows, assume that $\{\zeta _{t},t\geq 0\}$ has been chosen to be
this version. Given this version $\zeta $, Proposition \ref%
{prop:filterEqns:kallstrie} suggests the following definition:

\begin{definition}
\label{def:filterEqns:ucd} Define the \emph{unnormalised conditional
distribution} of $X$ to be the measure-valued process $\pi =\{\pi _{t},\
t\geq 0\}\ $given by $\pi _{t}=\zeta _{t}\varsigma _{t}$ for any $t\geq 0.$
\end{definition}

\paragraph{Notation.}
 We shall denote by $\mathcal{M}_F(\mathbb{R}^d)$ the set of
finite measures on $\mathbb{R}^d$, which we equip with the topology of weak
convergence (i.e. $\mu_n\to\mu$ if $\langle\mu_n,\varphi\rangle\to\langle%
\mu,\varphi\rangle$, for all $\varphi\in C_b(\mathbb{R}^d)$), the set of continuous bounded functions on $\mathbb{R}^d$.

\begin{lemma}
\label{le:filterEqns:rhoCadlag} Under assumption {\bf M}, the process $\{\pi _{t},\ t\geq 0\}$ is an 
$\mathcal{M}_F(\mathbb{R}^d)$--valued c\`{a}dl\`{a}g and $%
\mathcal{P}(\mathcal{Y}_{t})$-measurable process. Furthermore, for any $t\geq 0$, $%
\varphi\in b\mathcal{B}(\mathbb{R}^d)$, 
\begin{equation}
\pi _{t}(\varphi )=\tilde{\mathbb{E}}\left[ \tilde{Z}_{t}\varphi (X_{t})\mid 
\mathcal{Y}\right] \quad \quad \tilde{\mathbb{P}}(\mathbb{P})\text{-a.s.}
\label{eq:filterEqns:altRho}
\end{equation}
\end{lemma}

\begin{proof}
Both $\varsigma _{t}(\varphi )$ and $\zeta _{t}$ are $\mathcal{P}(\mathcal{Y}_{t}$)--measurable. 
By construction $\{\zeta _{t},\ t\geq 0\}$ is also c\`{a}dl\`{a}g.
Moreover, there exists a suitable version of the process $\varsigma
=\{\varsigma _{t},\ t\geq 0\}$, so that $\varsigma _{t}$ is 
$\mathcal{P}(\mathcal{Y}_{t})$-measurable probability measure-valued process for which (\ref%
{nfp}) holds almost surely, see Theorem 2.24 in \cite{bc}. 
In addition,
since $\mathcal{Y}_{t}$ is right-continuous, it follows that $\varsigma $
has a c\`adl\`ag version (see Corollary 2.26 in \cite{bc}). In the
following, we take $\varsigma $ to be this version. Moreover, for any
continuous bounded function $\varphi$, $\varsigma_{t}(\varphi )$ is the
optional projection of $\varphi (X_{t})$ with respect to the filtration $%
\mathcal{Y}_{t}$. Finally, since $\{\varsigma _{t},\ t\geq 0\}$ is c\`{a}dl%
\`{a}g and $\mathcal{P}(\mathcal{Y}_t)$ measurable, it follows that the process $\{
\pi_{t},t\geq 0\}$ is c\`{a}dl\`{a}g and $\mathcal{P}(\mathcal{Y}_t)$-measurable.

For the second part, from Proposition \ref{prop:filterEqns:p4} and
Proposition \ref{prop:filterEqns:kallstrie} it follows that 
\begin{equation*}
\pi _{t}(\varphi )=\varsigma _{t}(\varphi )\tilde{\mathbb{E}}[\tilde{Z}%
_{t}\mid \mathcal{Y}]=\tilde{\mathbb{E}}[\tilde{Z}_{t}\varphi (X_{t})\mid 
\mathcal{Y}]\quad \quad \tilde{\mathbb{P}}\text{-a.s.},
\end{equation*}%
From (\ref{eq:filterEqns:zeta}), $\tilde{\mathbb{E}}[\tilde{Z}_{t}\mid 
\mathcal{Y}_{t}]=\zeta _{t}$ a.s. from which the result follows.
\end{proof}

Definition \ref{def:filterEqns:ucd} gives us the following immediate
corollary:

\begin{corollary}
Under assumption {\bf M}, for every $\varphi \in b\mathcal{B}(\mathbb{R}^d)$, 
\begin{equation}
\varsigma _{t}(\varphi )=\frac{\pi _{t}(\varphi )}{\pi _{t}(\mathbf{1})}%
\quad \quad \forall t\in \lbrack 0,\infty )\quad \quad \tilde{\mathbb{P}}(%
\mathbb{P})\text{-a.s.}  \label{eq:filterEqns:corKS}
\end{equation}
\end{corollary}

\begin{remark}

The fact that the process $\varsigma$ has c\`adl\`ag paths is an application
of the properties of the optional projection of a stochastic process, see \cite{bc} for
further details. Whilst it is true that the optional projection of a process
with c\`adl\`ag paths has c\`adl\`ag paths, it is not, in general, true that
the optional projection of a continuous process is a continuous process. A
counterexample is the Azema martingale, see Theorem 61, pp 180-182, \cite%
{protter2005stochastic}.\footnote{%
We thank Martin Clark for pointing out this example to us.}

Continuity of both $\varsigma $ and $\pi $ can be ensured if additional
constraints are imposed. For example, if we have that 
\begin{equation}
\tilde{\mathbb{E}}[\sup_{t\in \lbrack 0,T]}\tilde{Z}_{t}]<\infty ,
\label{supremum}
\end{equation}%
for arbitrary $T>0$, then the process $\zeta $ is continuous by the
(conditional) dominated convergence theorem. Similarly, the measure valued
process $\{\pi _{t},\ t\geq 0\}$ also has continuous paths in $\mathcal{M}(%
\mathbb{R}^{d})$, the space of finite measures endowed with the weak
topology. In turn this implies that also the process $\{\varsigma _{t},\
t\geq 0\}$ has continuous paths in the same topology by (\ref%
{eq:filterEqns:corKS}).
\end{remark}

The Kallianpur--Striebel formula explains the usage of the term \emph{%
unnormalised} in the definition of $\pi _{t}$ as the denominator $\pi _{t}(%
\mathbf{1})$ can be viewed as the normalizing factor. Below $k^{+}$ stands
for the Moore Penrose pseudo-inverse of the matrix $k$ (see the Appendix below for its definition). Since
$(\omega,t)\mapsto k(t,Y_t(\omega))$ is $\mathcal{P}(\mathcal{Y}_t)$-measurable, and from Lemma \ref{MP-measurable} $k\mapsto k^+$ is measurable, we have that $k^+(t,Y_t)$ is 
$\mathcal{P}(\mathcal{Y}_t)$ measurable. Finally  $h_{2}^{\top }$
stands for the transpose of $h_{2}$ in other words the row vector
corresponding to $h_{2}$

\begin{lemma}
\label{integrarho} For all $t\geq 0,$ we have that 
\begin{equation}
	\int_{0}^{t}\left( \left\vert \left\vert \pi _{s}(\bar{g})k^{+}k\left(
	s,Y_{s}\right) \right\vert \right\vert ^{2}+\left\vert \left\vert \pi
	_{s}(h_{2}^{\top })k^{+}k\left( s,Y_{s}\right) \right\vert \right\vert
	^{2}\,\right) \mathrm{d}s<\infty ,\ \ \ \tilde{\mathbb{P}}\mathrm{-a.s.}
	\label{adcondonh}
\end{equation}%
As a consequence, the stochastic integrals
\begin{align*}
t\mapsto& \int_{0}^{t}\pi _{s}(\nabla \varphi _{s}\bar{g}+\varphi_{s}h_{2}^{\top })k^{+}\left( s,Y_{s}\right) (\mathrm{d}Y_{s}-h_1(s,Y_s)\mathrm{d}s),\quad \text{and}\\
t\mapsto& \int_{0}^{t}\left( \varsigma _{s}(\nabla \varphi _{s}\bar{g}+\varphi
_{s}h_{2}^{\top })-\varsigma _{s}(\varphi )\varsigma _{s}(h_{2}^{\top})\right) k^{+}\left( s,Y_{s}\right) \left( \mathrm{d}Y_{s}-\varsigma_{s}(h)\,\mathrm{d}s\right)
\end{align*}
are well defined for any function $\varphi \in C_{b}^{1,2}([0,\infty )\times 
\mathbb{R}^{d}))$ and are local semi-martingales with almost surely
continuous paths.
\end{lemma}

\begin{proof}
We only treat the first term in (\ref{adcondonh}). The second can be dealt
with in the same way. We first note that 
\begin{equation*}
\pi _{s}(\bar{g})k^{+}k\left( s,Y_{s}\right) =\varsigma _{s}(\bar{g})\pi
_{s}(1)k^{+}k\left( s,Y_{s}\right) ,
\end{equation*}%
and for any $0\leq s\leq t$, 
\begin{equation*}
\left\Vert \varsigma _{s}(\bar{g})\pi _{s}(1)k^{+}k\left( s,Y_{s}\right)
\right\Vert \leq \sup_{0\leq s\leq t}\pi _{s}(1)\left\Vert k^{+}k\left(
s,Y_{s}\right) \right\Vert \times \left\Vert \varsigma _{s}(\bar{g}%
)\right\Vert 
\end{equation*}%
Since $k^{+}k$ is locally bounded and the supremum on the interval $[0,t]$ of
the process $\pi _{s}(1)$ ($\pi _{s}$ has c\'{a}dl\'{a}g paths) and of the
process $\left\Vert Y_{s}\right\Vert $ ($Y$ has continuous paths) are finite 
$\tilde{P}$-a.s., clearly from \eqref{k+k} below,
\begin{equation*}
\sup_{0\leq s\leq t}\pi _{s}(1)^{2}\left\Vert k^{+}k\left( s,Y_{s}\right)
\right\Vert ^{2}<\infty ,\ \ \ {\tilde{P}}-\text{ a.s. }
\end{equation*}%
So it suffices to show that 
\begin{equation*}
\int_{0}^{t}\left\Vert \varsigma _{s}(\bar{g})\right\Vert ^{2}ds<\infty 
\text{ a.s. }
\end{equation*}%
This, in turn, is a consequence of the fact that for any entry $\bar{g}_{i,j}
$ of the matrix $\bar{g}$, using Jensen's inequality, 
\begin{equation*}
\begin{aligned} \mathbb{E} \int_0^t\left|\varsigma_s\left(\bar{g}_{i,
j}\right)\right|^2 d s & =\mathbb{E} \int_0^t \mid
\mathbb{E}\left[\bar{g}_{i, j}\left|\mathcal{Y}_s\right|^2 |d s\right. \\ &
\leq \mathbb{E} \int_0^t \mathbb{E}\left[\left|\bar{g}_{i, j}\right|^2 \mid
\mathcal{Y}_s\right] d s \\ & =\mathbb{E} \int_0^t\left|\bar{g}_{i,
j}\left(s, X_s, Y_s\right)\right|^2 d s \\ & <\infty, \end{aligned}
\end{equation*}%
where the last inequality follows from the fact that all coefficients have
at most linear growth and (\ref{squareint}). 

We next consider the two stochastic integrals. For the first one, we note that it can be rewritten as
\begin{align*}
\int_{0}^{t}\pi _{s}(\nabla \varphi _{s}\bar{g}+\varphi_{s}h_{2}^{\top })k^{+}k\left( s,Y_{s}\right) \mathrm{d}\tilde{W}_s,
\end{align*}
and the second one equals
\begin{align*}
\int_{0}^{t}\left( \varsigma _{s}(\nabla \varphi _{s}\bar{g}+\varphi
_{s}h_{2}^{\top })-\varsigma _{s}(\varphi )\varsigma _{s}(h_{2}^{\top})\right) k^{+}k\left( s,Y_{s}\right)\left(\mathrm{d}\tilde{W}_s-\varsigma _{s}(h_2)\mathrm{d}s\right),
\end{align*}
so that the result follows from the first part of the proof, combined with the fact that $\|k^{+}k\left( s,Y_{s}\right)\|\le C$, see \eqref{k+k} below.
\end{proof}

We are now in a position to establish the evolution equations for the
processes $\pi $ and $\varsigma $ in this set-up. Note that all stochastic
integrals in the equations for $\pi $ are well defined as per the above
Lemma. Similarly all the deterministic integrals are well defined as the
integrands are locally bounded. Similar arguments apply to the integrals in
the equation for $\varsigma $.

\begin{theorem}
\label{th:ex} Under assumption  \textbf{E}, the process $\pi _{t}$ satisfies
the following evolution equation 
\begin{eqnarray}
	\pi _{t}(\varphi _{t}) &=&\pi _{0}(\varphi _{0})+\int_{0}^{t}\pi _{s}\left(
	\partial _{s}\varphi _{s}+A_{s}\varphi _{s}\right) \,\mathrm{d}s  \notag \\
	&&+\int_{0}^{t}\pi _{s}(\nabla \varphi _{s}\bar{g}+\varphi _{s}h_{2}^{\top
	})k^{+}\left( s,Y_{s}\right) \left( \mathrm{d}Y_{s}-h_{1}\left(
	s,Y_{s}\right) \mathrm{d}s\right) ,\ \quad \tilde{\mathbb{P}}\text{-a.s.},\
	\ \forall t\geq 0  \label{eq:filterEqns:zakaifunctint}
\end{eqnarray}%
for any function $\varphi \in C_{b}^{1,2}([0,\infty )\times \mathbb{R}%
^{d})$\footnote{The set $C_{b}^{1,2}([0,\infty )\times \mathbb{R}%
	^{d})$ is the set of functions $\varphi : [0,\infty )\times \mathbb{R}
	^{d} \rightarrow \mathbb{R}$ that are once differentiable in the first variable and twice differentiable in the second variable and have all derivatives bounded.}. Moreover the conditional distribution $\varsigma _{t}$ satisfies
the following evolution equation 
\begin{align}
\varsigma _{t}(\varphi _{t})={}& \varsigma _{0}(\varphi
_{0})+\int_{0}^{t}\varsigma _{s}(\partial _{s}\varphi _{s}+A_{s}\varphi
_{s})\,\mathrm{d}s  \notag \\
& +\int_{0}^{t}\left( \varsigma _{s}(\nabla \varphi _{s}\bar{g}+\varphi
_{s}h_{2}^{\top })-\varsigma _{s}(\varphi )\varsigma _{s}(h_{2}^{\top
})\right) k^{+}\left( s,Y_{s}\right) \left( \mathrm{d}Y_{s}-\varsigma
_{s}(h)\,\mathrm{d}s\right)   \label{eq:filterEqns:ks}
\end{align}
for any function $\varphi \in C_{b}^{1,2}([0,\infty )\times \mathbb{R}%
^{d})$.
\end{theorem}

\begin{remark}
Assumption \textbf{E} includes the degenerate case $k=0$. In this case, the
observation process satisfies the evolution equation 
\begin{equation}
Y_{t} =Y_{0}+\int_{0}^{t}h_1(s,Y_s)\,\mathrm{d}s.  \label{Yfork=0}
\end{equation}
It follows that $Y$ being deterministic, it is independent of $X$. Hence
both $\varsigma$ and $\pi$ coincide with the law of the signal $X$. In
particular, since $k^+=0^+=0$, equation (\ref{eq:filterEqns:zakaifunctint})
degenerates to

\begin{equation}
\pi _{t}(\varphi _{t})=\pi _{0}(\varphi _{0})+\int_{0}^{t}\pi _{s}\left(
\partial _{s}\varphi _{s}+A_{s}\varphi _{s}\right) \,\mathrm{d}s
\label{eq:deg_filterEqns:zakaifunctint}
\end{equation}%
and so does equation (\ref{eq:filterEqns:ks}).
\end{remark}

\noindent \emph{\ Proof of Theorem \ref{th:ex}.}

We can re-write the equation satisfied by the process $(X,Y,Z)$ as being
driven by the pair of processes $(V,\tilde{W})$ 
\begin{eqnarray*}
X_{t} &=&\!X_{0}\!+\!\int_{0}^{t}[f(s,X_{s},Y_{s})-\bar{g}%
h_{2}(s,X_{s},Y_{s})]\mathrm{d}s\!+\!\int_{0}^{t}\!\!g(s,X_{s},Y_{s})\,%
\mathrm{d}V_{s}+\int_{0}^{t}\!\!\bar{g}(s,X_{s},Y_{s})\mathrm{d}\tilde{W}_{s}
\\
Y_{t} &=&Y_{0}+\int_{0}^{t}h_{1}(s,Y_{s})\,\mathrm{d}s+\int_{0}^{t}k(s,Y_{s})%
\mathrm{d}{\tilde{W}}_{s} \\
\tilde{Z}_{t} &=&1+\int_{0}^{t}\tilde{Z}_{s}\left(
h_{2}(s,X_{s},Y_{s})\right) ^{\top }\,\mathrm{d}\tilde{W}_{s}
\end{eqnarray*}%
To ensure the integrability of the terms appearing in the following
computations, we first approximate $\tilde{Z}_{t}$ with $\tilde{Z}%
_{t}^{\varepsilon }$ given by 
\begin{equation*}
\tilde{Z}_{t}^{\varepsilon }=\frac{\tilde{Z}_{t}}{1+\varepsilon \tilde{Z}_{t}%
}=\frac{1}{\varepsilon }\frac{\varepsilon \tilde{Z}_{t}}{1+\varepsilon 
\tilde{Z}_{t}}=\frac{1}{\varepsilon }\left( 1-\frac{1}{1+\varepsilon \tilde{Z%
}_{t}}\right) .
\end{equation*}
By It\^{o}'s formula, we deduce that\footnote{%
In the following $\nabla \varphi \,$will denote the row vector $\left(
\partial _{1}\varphi ,...,\partial _{d}\varphi \right) $.} 
\begin{eqnarray*}
\mathrm{d}\varphi \left( t,X_{t}\right) &=&[(\partial _{t}+A_{t})\varphi
](t,X_{t})\mathrm{d}t-\nabla\varphi \left( t,X_{t}\right) \bar{g}
h_{2}(t,X_{t},Y_{t})]\mathrm{d}t \\
&&+(\nabla \varphi g)(t,X_{t})\,\mathrm{d}V_{t}+\nabla \varphi \bar{g}
\,(t,X_{t})\,\mathrm{d}\tilde{W}_{t} \\
\mathrm{d}\tilde{Z}_{t}^{\varepsilon } &=&\tilde{Z}_{t}(1+\varepsilon \tilde{
Z}_{t})^{-2}\left( h_{2}(t,X_{t},Y_{t})\right) ^{\top }\,\mathrm{d}\tilde{W}%
_{t} \\
&&-\varepsilon (1+\varepsilon \tilde{Z}_{t})^{-3}\tilde{Z}_{t}^{2}\left(
h_{2}(t,X_{t},Y_{t})\right) ^{\top }h_{2}(t,X_{t},Y_{t})\,\mathrm{d}t
\end{eqnarray*}%
Therefore%
\begin{eqnarray}
\mathrm{d}\varphi \left( t,X_{t}\right) \tilde{Z}_{t}^{\varepsilon } &=&%
\tilde{Z}_{t}^{\varepsilon }[(\partial _{t}+A_{t})\varphi
](t,X_{t})-\nabla\varphi \left( t,X_{t}\right) \bar{g}h_{2}(t,X_{t},Y_{t})]%
\mathrm{d}t  \notag \\
&&+\tilde{Z}_{t}^{\varepsilon }\left( (\nabla \varphi g)(t,X_{t})\,\mathrm{d}%
V_{t}+\nabla \varphi \bar{g}\,(t,X_{t})\,\mathrm{d}\tilde{W}_{t}\right) 
\notag \\
&&+\varphi \left( t,X_{t}\right) \tilde{Z}_{t}(1+\varepsilon \tilde{Z}%
_{t})^{-2}\left( h_{2}(t,X_{t},Y_{t})\right) ^{\top }\,\mathrm{d}\tilde{W}%
_{t}  \notag \\
&&+\varphi \left( t,X_{t}\right) \left( -\varepsilon (1+\varepsilon \tilde{Z}%
_{t})^{-3}\tilde{Z}_{t}^{2}\left( h_{2}(t,X_{t},Y_{t})\right) ^{\top
}h_{2}(t,X_{t},Y_{t})\,\mathrm{d}t\right)  \notag \\
&&+\nabla \varphi \bar{g}\,(t,X_{t})h_{2}(t,X_{t},Y_{t})\tilde{Z}%
_{t}(1+\varepsilon \tilde{Z}_{t})^{-2}\mathrm{d}t  \notag \\
&=&\tilde{Z}_{t}^{\varepsilon }[(\partial _{t}+A_{t})\varphi ](t,X_{t})%
\mathrm{d}t-\nabla\varphi \left( t,X_{t}\right) \bar{g}h_{2}(t,X_{t},Y_{t})]%
\mathrm{d}t  \notag \\
&&+\tilde{Z}_{t}^{\varepsilon }\nabla \varphi \bar{g}\,(t,X_{t})\,k^{+}k%
\left( t,Y_{t}\right) \,\mathrm{d}\tilde{W}_{t}  \notag \\
&&+\tilde{Z}_{t}^{\varepsilon }\nabla \varphi \bar{g}\,(t,X_{t})\,(I-k^{+}k%
\left( t,Y_{t}\right) )\,\mathrm{d}\tilde{W}_{t}  \notag \\
&&+\tilde{Z}_{t}^{\varepsilon }(\nabla \varphi g)(t,X_{t})\,\mathrm{d}V_{t} 
\notag \\
&&+\varphi \left( t,X_{t}\right) \tilde{Z}_{t}^{\varepsilon }(1+\varepsilon 
\tilde{Z}_{t})^{-1}\left( h_{2}(t,X_{t},Y_{t})\right) ^{\top }\,k^{+}k\left(
t,Y_{t}\right) \mathrm{d}\tilde{W}_{t}  \notag \\
&&+\varphi \left( t,X_{t}\right) \tilde{Z}_{t}^{\varepsilon }(1+\varepsilon 
\tilde{Z}_{t})^{-1}\left( h_{2}(t,X_{t},Y_{t})\right) ^{\top
}(I-k^{+}k\left( t,Y_{t}\right) )\mathrm{d}\tilde{W}_{t}  \notag \\
&&+\varphi \left( t,X_{t}\right) \left( -\varepsilon (1+\varepsilon \tilde{Z}%
_{t})^{-1}\left( \tilde{Z}_{t}^{\varepsilon }\right) ^{2}\left(
h_{2}(t,X_{t},Y_{t})\right) ^{\top }h_{2}(t,X_{t},Y_{t})\,\mathrm{d}t\right)
\notag \\
&&+\nabla \varphi \bar{g}\,(t,X_{t})h_{2}(t,X_{t},Y_{t})\tilde{Z}%
_{t}^{\varepsilon }(1+\varepsilon \tilde{Z}_{t})^{-2}\mathrm{d}t
\label{longform}
\end{eqnarray}

We next take the conditional expectation $\tilde{\mathbb{E}}%
(\cdot |\mathcal{Y})$ in this identity, as all the terms in (\ref{longform})
are square integrable over $\left[ 0,T\right] \times \Omega .$ To show this
we use repeatedly the fact that 
\[
| \tilde{Z}_{t}^{\varepsilon }|=| \tilde{Z}%
_{t}(1+\varepsilon \tilde{Z}_{t})^{-1}| \leq \varepsilon ^{-1}, \ \  |
(1+\varepsilon \tilde{Z}_{t})^{-1} | \leq 1, 
\]
that $\varphi \in
C_{b}^{1,2}([0,\infty )\times \mathbb{R}^{d})),$ the linear growth of $f,g,%
\bar{g},h_{2}$ and $k$, and the square integrability of the processes $X$ and 
$Y$.

We next take the conditional expectation $\tilde{\mathbb{E}}(\cdot |\mathcal{%
Y})$ in this identity, as all the terms in (\ref{longform}) are integrable
over $\left[ 0,T\right] \times \Omega .$  The conditional expectation  operator
$\tilde{\mathbb{E}}(\cdot |\mathcal{Y})$ commutes with the $\mathrm{d}t$ and
the $\mathrm{d}\tilde{W}_{t}$ integrations, while $\tilde{\mathbb{E}}(\cdot |%
\mathcal{Y})$ of a stochastic integral w.r.t. $dV_{S}$ and to $%
[I-k^{+}k(s,Y_{s})]\,\mathrm{d}\tilde{W}_{s}$ is zero. The last claim is the
content of the Lemma \ref{lem:cond}. We deduce that
\begin{eqnarray}
\tilde{\mathbb{E}}[\tilde{Z}_{t}^{\varepsilon }\varphi (t,X_{t})\mid 
\mathcal{Y}] &=&{}\frac{\pi _{0}(\varphi )}{1+\varepsilon }+\int_{0}^{t}%
\tilde{\mathbb{E}}\left[ \tilde{Z}_{s}^{\varepsilon }[(\partial
_{s}+A_{s})\varphi ](s,X_{s})-\nabla \varphi \left( s,X_{s}\right) \bar{g}%
h_{2}(s,X_{s},Y_{s})]\mid \mathcal{Y}\right] \,\mathrm{d}s  \notag \\
&&+\int_{0}^{t}\tilde{\mathbb{E}}[\tilde{Z}_{s}^{\varepsilon }\nabla \varphi 
\bar{g}\,(s,X_{s})\,|\mathcal{Y}]\,\,k^{+}k\left( s,Y_{s}\right) \mathrm{d}%
\tilde{W}_{s}  \notag \\
&&+\int_{0}^{t}\tilde{\mathbb{E}}[\varphi \left( s,X_{s}\right) \tilde{Z}%
_{s}^{\varepsilon }(1+\varepsilon \tilde{Z}_{s})^{-1}\left(
h_{2}(s,X_{s},Y_{s})\right) ^{\top }\,\,|\mathcal{Y}]\,\,k^{+}k\left(
s,Y_{s}\right) \mathrm{d}\tilde{W}_{s}  \notag \\
&&+\int_{0}^{t}\tilde{\mathbb{E}}\left[ \varphi \left( s,X_{s}\right) \left(
-\varepsilon (1+\varepsilon \tilde{Z}_{s})^{-1}\left( \tilde{Z}%
_{s}^{\varepsilon }\right) ^{2}\left( h_{2}(s,X_{s},Y_{s})\right) ^{\top
}h_{2}(s,X_{s},Y_{s})\,\right) |\ \mathcal{Y}\right] \,\,\mathrm{d}s  \notag
\\
&&+\int_{0}^{t}\tilde{\mathbb{E}}[\nabla \varphi \bar{g}%
\,(s,X_{s})h_{2}(s,X_{s},Y_{s})\tilde{Z}_{s}^{\varepsilon }(1+\varepsilon 
\tilde{Z}_{s})^{-2}\,|\mathcal{Y}]\,\mathrm{d}s  \label{lf2}
\end{eqnarray}

Using Proposition \ref{prop:filterEqns:p4} we deduce from (\ref{lf2}) by
taking the limit as that $\varepsilon $ tends to 0 that 
\begin{eqnarray}
\pi _{t}(\varphi _{t}) &=&\pi _{0}(\varphi _{0})+\int_{0}^{t}\pi _{s}\left(
\partial _{s}\varphi _{s}+A_{s}\varphi _{s}-\nabla \varphi \bar{g}%
h_{2}+\nabla \varphi \bar{g}h_{2}\right) \,\mathrm{d}s+  \notag \\
&&+\int_{0}^{t}\pi _{s}(\nabla \varphi _{s}\bar{g}+\varphi _{s}h_{2}^{\top
})k^{+}\left( s,Y_{s}\right) \left( \mathrm{d}Y_{s}-h_{1}\left(
s,Y_{s}\right) \mathrm{d}s\right) ,\ \quad \tilde{\mathbb{P}}\text{-a.s.},\
\ \forall t\geq 0  \label{lf4}
\end{eqnarray}

 In order to justify taking that the limit in (\ref{lf2}) gives (%
\ref{lf4}), we need to show that the integrands on the right hand side of (%
\ref{lf2}) are uniformly bounded in $\varepsilon $ by processes that are
integrable over the product space $\left[ 0,t\right] \times \Omega $. First
we have that 
\begin{eqnarray*}
&&\left\vert \tilde{Z}_{s}^{\varepsilon }[(\partial _{s}+A_{s})\varphi
](s,X_{s})-\nabla \varphi \left( s,X_{s}\right) \bar{g}h_{2}(s,X_{s},Y_{s})]
\right\vert \\
& &\leq \tilde{Z}\left\vert [(\partial _{s}+A_{s})\varphi ](s,X_{s})-\nabla
\varphi \left( s,X_{s}\right) \bar{g}h_{2}(s,X_{s},Y_{s})] \right\vert \\
&&\leq c_1\tilde{Z}\left( 1+\sum_{i=1}^{d}\left\vert f^{i}\left(
s,X_{s},Y_{s}\right) \right\vert +\frac{1}{2}\sum_{i,j=1}^{d}\left\vert
a^{ij}\left( s,X_{s},Y_{s}\right) \right\vert
+\sum_{i=1}^{d}\sum_{j=1}^{l^{\prime }}\left\vert \bar{g}^{ij}\left(
s,X_{s},Y_{s}\right) h_{2}^{j}\left( s,X_{s},Y_{s}\right) \right\vert \right)
\\
&&\leq c_2\tilde{Z}\left( 1+\left\vert \left\vert X_{s}\right\vert
\right\vert ^{2}+\left\vert \left\vert Y_{s}\right\vert \right\vert
^{2}\right) ,
\end{eqnarray*}%
where we used the fact that $\varphi \in C_{b}^{1,2}([0,\infty )\times 
\mathbb{R}^{d}))$, the linear growth of $f,g,\bar{g},h_{2}$ and $k$ and
denoted by $c_1=c_1\left( \varphi \right) =\left\vert \left\vert \partial
_{s}\varphi \right\vert \right\vert _{\infty }+\left\vert \left\vert \nabla
\varphi \right\vert \right\vert _{\infty }+\left\vert \left\vert \nabla
\nabla \varphi \right\vert \right\vert _{\infty }$ and $c_2$ is a constant
depending of $c_1$ and on the constant $K$ from Assumption E. We then
observe that 
\begin{equation*}
\int_{0}^{t}\tilde{\mathbb{E}}\left[ c\tilde{Z}\left( 1+\left\vert
\left\vert X_{s}\right\vert \right\vert ^{2}+\left\vert \left\vert
Y_{s}\right\vert \right\vert ^{2}\right) \right] \,\mathrm{d}s=c\int_{0}^{t} 
\mathbb{E}\left[ \left( 1+\left\vert \left\vert X_{s}\right\vert \right\vert
^{2}+\left\vert \left\vert Y_{s}\right\vert \right\vert ^{2}\right) \right]
\,\mathrm{d}s<\infty ,
\end{equation*}%
which justifies that the integrand in the first term on the right hand side
of (\ref{lf2}) is dominated by an integrable bound independent of $%
\varepsilon $. This justifies the convergence of the first term in (\ref{lf2}%
). 

A similar argument applies to the last two terms of (\ref{lf2}), using

\begin{eqnarray*}
\left\vert \varphi \left( s,X_{s}\right) \left( -\varepsilon (1+\varepsilon 
\tilde{Z}_{s})^{-1}\left( \tilde{Z}_{s}^{\varepsilon }\right) ^{2}\left(
h_{2}(s,X_{s},Y_{s})\right) ^{\top }h_{2}(s,X_{s},Y_{s})\,\right)
\right\vert &\leq &c\tilde{Z}\left( 1+\left\vert \left\vert X_{s}\right\vert
\right\vert ^{2}+\left\vert \left\vert Y_{s}\right\vert \right\vert
^{2}\right) \\
\left\vert \nabla \varphi \bar{g}\,(s,X_{s})h_{2}(s,X_{s},Y_{s})\tilde{Z}%
_{s}^{\varepsilon }(1+\varepsilon \tilde{Z}_{s})^{-1}\right\vert &\leq &c%
\tilde{Z}\left( 1+\left\vert \left\vert X_{s}\right\vert \right\vert
^{2}+\left\vert \left\vert Y_{s}\right\vert \right\vert ^{2}\right)
\end{eqnarray*}%
The convergence of the stochastic terms is harder. We combine them into a
single term and re-write it as 
\begin{equation}
M_{t}^{\varepsilon }=\int_{0}^{t}q_{s}^{\varepsilon }k^{+}\left(
s,Y_{s}\right) \left( \mathrm{d}Y_{s}-h_{1}\left( s,Y_{s}\right) \mathrm{d}%
s\right) =:\int_{0}^{t}q_{s}^{\varepsilon }k^{+}k\left( s,Y_{s}\right) 
\mathrm{d}\tilde{W}_{s},~~~t\geq 0,  \label{sumsto}
\end{equation}%
where 
\begin{equation*}
q_{s}^{\varepsilon }=\tilde{\mathbb{E}}[\tilde{Z}_{s}^{\varepsilon }\nabla
\varphi \bar{g}\,(s,X_{s})+\varphi \tilde{Z}_{s}^{\varepsilon
}(1+\varepsilon \tilde{Z}_{s})^{-1}\left( h_{2}(s,X_{s},Y_{s})\right) ^{\top
}\,|\mathcal{Y}],~~s\geq 0
\end{equation*}%
Observe that $M^{\varepsilon }$\ is a square integrable martingale, however
the intended limit%
\begin{eqnarray*}
M_{t}:= &&\int_{0}^{t}q_{s}k^{+}k\left( s,Y_{s}\right) \mathrm{d}\tilde{W}%
_{s}=\int_{0}^{t}\pi _{s}(\nabla \varphi _{s}\bar{g}+\varphi _{s}h_{2}^{\top
})k^{+}\left( s,Y_{s}\right) \left( \mathrm{d}Y_{s}-h_{1}\left(
s,Y_{s}\right) \mathrm{d}s\right) ,~~~,~~~t\geq 0, \\
q_{s}:= &&\pi _{s}(\nabla \varphi _{s}\bar{g}+\varphi _{s}h_{2}^{\top
}),~~~~s\geq 0
\end{eqnarray*}%
is only a local martingale for any function $\varphi \in
C_{b}^{1,2}([0,\infty )\times \mathbb{R}^{d}))$ with almost surely
continuous paths which  may not be  square integrable, (see Lemma \ref%
{integrarho} for details). To begin the convergence argument, observe that for any $0\le s\le t$ and a.s. 
\begin{equation*}
\lim_{\varepsilon \rightarrow \infty }\tilde{Z}_{s}^{\varepsilon }\nabla
\varphi \bar{g}\,(s,X_{s})+\varphi \tilde{Z}_{s}^{\varepsilon
}(1+\varepsilon \tilde{Z}_{s})^{-1}\left( h_{2}(s,X_{s},Y_{s})\right) ^{\top
}\,=\tilde{Z}_{s}\nabla \varphi \bar{g}\,(s,X_{s})+\varphi \tilde{Z}%
_{s}\left( h_{2}(s,X_{s},Y_{s})\right) ^{\top }\,
\end{equation*}%
and 
\begin{equation}
\left\vert \left\vert \tilde{Z}_{s}^{\varepsilon }\nabla \varphi \bar{g}%
\,(s,X_{s})+\varphi \tilde{Z}_{s}^{\varepsilon }(1+\varepsilon \tilde{Z}%
_{s})^{-1}\left( h_{2}(s,X_{s},Y_{s})\right) ^{\top }\right\vert \right\vert
\leq c\tilde{Z}_{s}\sqrt{\left( 1+\left\vert \left\vert X_{s}\right\vert
\right\vert ^{2}+\left\vert \left\vert Y_{s}\right\vert \right\vert
^{2}\right) }.  \label{integr1}
\end{equation}%
Therefore, since the term on the right hand side of (\ref{integr1}) is $\tilde{\mathbb{P}}$
integrable we deduce by the dominated convergence theorem for conditional expectation that 
\begin{equation*}
\lim_{\varepsilon \rightarrow \infty }q_{s}^{\varepsilon }=q_{s}
\end{equation*}%
$\tilde{P}$-almost surely and almost everywhere on the interval $\left[ 0,t%
\right] $. Also one deduces that 
\begin{equation}
\left\vert \left\vert q_{s}^{\varepsilon }\right\vert \right\vert \leq
c\varsigma _{s}(\left\vert \left\vert \bar{g}\right\vert \right\vert
+\left\vert \left\vert h_{2}\right\vert \right\vert )\pi _{s}(1)
\label{integr2}
\end{equation}%
and since the term on the right hand side of (\ref{integr2}) is a.s. \emph{bounded%
} on the interval $\left[ 0,t\right]$ (the argument is similar to that used
in Lemma \ref{integrarho}), we deduce that (we use again the fact that $%
\left\vert \left\vert k^{+}k\right\vert \right\vert \leq C$ ) 
\begin{equation*}
0\leq \lim_{\varepsilon \rightarrow \infty }\int_{0}^{t}\left\vert
\left\vert \left( q_{s}^{\varepsilon }-q_{s}\right) k^{+}k\left(
s,Y_{s}\right) \right\vert \right\vert ^{2}ds\leq \int_{0}^{t}\left\vert
\left\vert q_{s}^{\varepsilon }-q_{s}\right\vert \right\vert ^{2}ds=0
\end{equation*}%
$\tilde{P}$-almost surely. This, in turn, implies that (for example by using
Proposition B.41. in \cite{bc}) 
\begin{equation*}
\lim_{\varepsilon \rightarrow 0}\sup_{t\in \left[ 0,T\right] }\left\vert
M_{t}^{\varepsilon}-M_{t}\right\vert =\lim_{\varepsilon \rightarrow
0}\sup_{t\in \left[ 0,T\right] }\left\vert \int_{0}^{t}\left(
q_{s}^{\varepsilon }-q_{s}\right) k^{+}k\left( s,Y_{s}\right) \mathrm{d}%
\tilde{W}_{s}\right\vert =0
\end{equation*}%
in probability. The justification
of the identity (\ref{eq:filterEqns:zakaifunctint}) is now complete.

To deduce that the conditional distribution of the signal $\varsigma _{t}$
satisfies (\ref{eq:filterEqns:ks}), we first compute the evolution equation
for the reciprocal of the mass process $\frac{1}{\pi _{t}(1)}$ which is 
\begin{eqnarray*}
\frac{1}{\pi _{t}(1)} &=&\frac{1}{\pi _{0}(1)}-\int_{0}^{t}\frac{\varsigma
_{s}(h_{2}^{\top })}{\pi _{s}(1)}k^{+}k\left( s,Y_{s}\right) \mathrm{d}{%
\tilde{W}}_{s}+\int_{0}^{t}\frac{\varsigma _{s}(h_{2}^{\top })}{\pi _{s}(1)}%
k^{+}k\left( k^{+}k\right) ^{\top }\left( s,Y_{s}\right) \varsigma
_{s}(h_{2})\mathrm{d}s \\
&=&\frac{1}{\pi _{0}(1)}-\int_{0}^{t}\frac{\varsigma _{s}(h_{2}^{\top })}{%
\pi _{s}(1)}k^{+}k\left( s,Y_{s}\right) \mathrm{d}{\tilde{W}}%
_{s}+\int_{0}^{t}\frac{\varsigma _{s}(h_{2}^{\top })}{\pi _{s}(1)}%
k^{+}k\left( s,Y_{s}\right) \varsigma _{s}(h_{2})\mathrm{d}s.
\end{eqnarray*}%
To obtain the second identity we used that, see \eqref{moorepenrose},
\begin{equation*}
\left( k^{+}k\right) ^{\top }=k^{+}k,~~~~~kk^{+}k=k\,.
\end{equation*}
Finally we use It\^o's formula to deduce that 
\begin{eqnarray*}
\varsigma _{t}(\varphi _{t}) &=&\frac{\pi _{t}(\varphi _{t})}{\pi _{t}(1)} \\
&=&\frac{\pi _{0}(\varphi _{0})}{\pi _{0}(1)}+\int_{0}^{t}\frac{\pi
_{s}\left( \partial _{s}\varphi _{s}+A_{s}\varphi _{s}\right) }{\pi _{s}(1)}%
\,\mathrm{d}s+\int_{0}^{t}\frac{\pi _{s}(\nabla \varphi _{s}\bar{g}+\varphi
_{s}h_{2}^{\top })}{\pi _{s}(1)}k^{+}k\left( s,Y_{s}\right) \mathrm{d}{%
\tilde{W}}_{s} \\
&&-\int_{0}^{t}\frac{\pi _{s}(\varphi _{s})\varsigma _{s}(h_{2}^{\top })}{%
\pi _{s}(1)}k^{+}k\left( s,Y_{s}\right) \mathrm{d}{\tilde{W}}%
_{s}+\int_{0}^{t}\frac{\pi_s(\varphi_s)\varsigma _{s}(h_{2}^{\top })}{\pi
_{s}(1)}k^{+}k\left( s,Y_{s}\right) \varsigma _{s}(h_{2})\mathrm{d}s \\
&&-\int_{0}^{t}\varsigma _{s}(\nabla \varphi _{s}\bar{g}+\varphi
_{s}h_{2}^{\top })k^{+}k\left( s,Y_{s}\right) \varsigma _{s}(h_{2})\mathrm{d}%
s \\
&=&\varsigma _{0}(\varphi _{0})+\int_{0}^{t}\varsigma _{s}(\partial
_{s}\varphi _{s}+A_{s}\varphi _{s})\,\mathrm{d}s \\
&&+\int_{0}^{t}\left( \varsigma _{s}(\nabla \varphi _{s}\bar{g}+\varphi
_{s}h_{2}^{\top })-\varsigma _{s}(\varphi )\varsigma _{s}(h_{2}^{\top
})\right) k^{+}\left( s,Y_{s}\right) \left( \mathrm{d}Y_{s}-\varsigma
_{s}(h)\,\mathrm{d}s\right) .
\end{eqnarray*}%
\qed

\begin{lemma}
\label{lem:cond} For any progressively measurable process $\kappa $ assumed
to be an $m$--dimensional row vector satisfying for any $t>0$ $\tilde{%
\mathbb{E}}\int_{0}^{t}\Vert \kappa _{s}\Vert ^{2}ds<\infty $, 
\begin{equation*}
\tilde{\mathbb{E}}\left[\int_{0}^{t}\kappa _{s}[I-k^{+}k(s,Y_{s})]\,\mathrm{%
d}\tilde{W}_{s}\Big|\mathcal{Y}\right] =0\,,\forall t>0\,.
\end{equation*}
\end{lemma}

\begin{proof}
Let $N=\left\{ N_{t},\ t\geq 0\right\} $ be the following process 
\begin{equation}
N_{t}=\int_{0}^{t}k(s,Y_{s})\,\mathrm{d}\tilde{W}_{s},~~~t\geq 0.
\label{ntdef}
\end{equation}%
We have
\begin{equation}
Y_{t}=Y_{0}+\int_{0}^{t}h_{1}(s,Y_{s})\,\mathrm{d}s+N_{t},\ t\ge0\,.  \label{NT}
\end{equation}%
Denote by $\mathcal{N}=\left\{ \mathcal{N}_{t},t\geq 0\right\} $ the
filtration generated by the process $N$.\ We prove that $\mathcal{Y}_{0}\vee 
\mathcal{N}_{t}=\mathcal{Y}_{t}$. To do this we deduce from (\ref{NT}) that $%
\mathcal{N}_{t}\subset $ $\mathcal{Y}_{t}$ and also that $\mathcal{Y}%
_{0}\vee \mathcal{N}_{t}$\ $\subseteq \mathcal{Y}_{t}$. Next, if we denote
by $\vartheta =\left\{ \vartheta _{t},\ t\geq 0\right\} $ the process $%
\vartheta _{t}=Y_{t}-N_{t}$, we get that 
\begin{equation*}
\vartheta _{t}=Y_{0}+\int_{0}^{t}h_{1}(s,\vartheta _{s}+N_{s})\,\mathrm{d}s.
\end{equation*}%
and we deduce from the above that the process $\vartheta $, which is the
unique solution of an ODE with initial condition $\vartheta _{0}=Y_{0}$, and
whose coefficient depends upon the process $N$ (note that $h_{1}$ is a
Lipschitz function in the spatial variable), is  
$
\mathcal{P}(\mathcal{Y}_{0}\vee \mathcal{N}_{t})$ measurable. From this we deduce that $Y_{t}$ is
measurable with respect to $\mathcal{Y}_{0}\vee \mathcal{N}_{t}$ and
therefore that $\mathcal{Y}_{t}\subseteq \mathcal{Y}_{0}\vee \mathcal{N}_{t}$%
, which in turn, implies that $\mathcal{Y}_{0}\vee \mathcal{N}_{t}=\mathcal{Y%
}_{t}$.

In the following we will use an argument first used by Krylov and Rozovsky
in \cite{KR}. Let us denote by $S_{t}$ the following set of uniformly
bounded test random variables: 
\begin{equation}
S_{t}=\left\{ \theta _{t}=\exp \!\left( i\bar{r}^{\top
}Y_{0}+i\int_{0}^{t}r_{s}^{\top }\,\mathrm{d}N_{s}+\frac{1}{2}%
\int_{0}^{t}\Vert k^{T}(s,Y_{s})r_{s}\Vert ^{2}\,\mathrm{d}s\right), \ r\in
L^{\infty }\left( [0,t],\mathbb{R}^{d^{\prime }}\right) ,~\bar{r}\in \mathbb{%
R}^{d^{\prime }}\right\} .  \label{class}
\end{equation}%
Then $S_{t}$ is a total set in $L^{1}(\Omega ,\mathcal{Y}_{0}\vee \mathcal{N}%
_{t},\tilde{\mathbb{P}})\equiv L^{1}(\Omega ,\mathcal{Y}_{t},\tilde{\mathbb{P%
}})$. That is, if $a\in L^{1}(\Omega ,\mathcal{Y}_{t},\tilde{\mathbb{P}})$
and $\tilde{\mathbb{E}}[a\theta _{t}]=0$, for all $\theta _{t}\in S_{t}$,
then $a=0$ $\tilde{\mathbb{P}}$-a.s. A proof of this result is similar with
that of Lemma B.39, pp 355 in  \cite{bc}. In addition, if $%
\theta _{t}\in S_{t}$, then 
\begin{equation*}
\theta _{t}=\exp \!\left( i\bar{r}^{\top }Y_{0}\right) +i\int_{0}^{t}\theta
_{s}r_{s}^{\top }\,\mathrm{d}N_{s}=\exp \!\left( i\bar{r}^{\top
}Y_{0}\right) +i\int_{0}^{t}\theta _{s}\,r_{s}^{T}k(s,Y_{s})\,\mathrm{d}\tilde{W}_{s}.
\end{equation*}
In view of this result, all we need to prove is that for any $r\in
L^{\infty}\left( [0,t],\mathbb{R}^{d^{\prime }}\right)$, $\bar{r}\in \mathbb{R}^{d^{\prime }}$, 
\begin{equation*}
\tilde{\mathbb{E}}\left( \theta _{t}\int_{0}^{t}\kappa
_{s}[I-k^{+}k(s,Y_{s})]\,\mathrm{d}\tilde{W}_{s}\right) =0\,.
\end{equation*}%
Since $k[I-k^{+}k]^{T}=0$, as a consequence of $(k^{+}k)^{T}=k^{+}k$ and $kk^{+}k=k$ (see \eqref{moorepenrose}), we have 
\begin{align*}
\tilde{\mathbb{E}}\left( \theta _{t}\int_{0}^{t}\kappa
_{s}[I-k^{+}k(s,Y_{s})]\,\mathrm{d}\tilde{W}_{s}\right) & =0+i\,\tilde{\mathbb{E}}\int_{0}^{t}\theta
_{s}r_{s}^{T}k(s,Y_{s})[I-k^{+}k(s,Y_{s})]^{T}\kappa _{s}^{T}ds \\
& =0.
\end{align*}
\end{proof}

\begin{remark}
In Lemma B.39 in \cite{bc}, the process  $Y$ is chosen to be a Brownian
motion starting from 0. In particular $Y_{0}=0$. \ However, the argument
does not use the property that $Y$ is a Brownian motion and can be extended
to any martingale. The additional term $i\bar{r}^{\top }Y_{0}$ in the
argument of $\theta _{t}$ takes care of the non-zero initial condition.
\end{remark}

The following lemma is an immediate consequence of Theorem 6 in Chapter V of \cite{protter2005stochastic}.

\begin{lemma}
\label{lemma:mass} Let $\varsigma $ be a solution of the evolution equation (%
\ref{eq:filterEqns:ks}) and consider the following stochastic differential
equation%
\begin{eqnarray}
dj_{t}^{\varsigma } &=&j_{t}^{\varsigma }\varsigma _{t}(h_{2}^{\top
})k^{+}\left( t,Y_{t}\right) dN_{t},  \notag \\
&=&j_{t}^{\varsigma }\varsigma _{t}(h_{2}^{\top })k^{+}k(t,Y_{t})\,\mathrm{d}%
\tilde{W}_{t}~~~j_{0}^{\varsigma }=1.  \label{abstracttotalmass}
\end{eqnarray}%
Then equation (\ref{abstracttotalmass}) has a unique solution $j^{\varsigma
}=\left\{ j_{t}^{\varsigma },t\geq 0\right\} $ for any solution $\varsigma $
of the evolution equation (\ref{eq:filterEqns:ks}).
\end{lemma}

 We next prove.
\begin{lemma}
\label{lemma:prob} Let $\varsigma $ be a solution of the evolution equation (%
\ref{eq:filterEqns:ks}) starting from $\varsigma _{0}$ which is a
probability measure. Then $\varsigma _{t}$ is a probability measure valued
process, P-almost surely. Moreover if $j^{\varsigma }$ is the corresponding
solution of equation (\ref{abstracttotalmass}), then $j_{t}^{\varsigma
}\varsigma _{t}$ is a solution of the evolution equation (\ref%
{eq:filterEqns:zakaifunctint}).
\end{lemma}

\begin{proof}
From (\ref{eq:filterEqns:ks}) we deduce that 
\begin{equation}
\varsigma _{t}(1)-1={}-\int_{0}^{t}\left( \varsigma _{s}(1)-1)\varsigma
_{s}(h_{2}^{\top })\right) k^{+}k(t,Y_{s})\,(\mathrm{d}\tilde{W}%
_{s}-\varsigma_s(h_2)ds),  \notag
\end{equation}
which implies that $a_{t}^{0}:=\varsigma _{t}(1)-1$ is a solution of the
linear equation 
\begin{equation}  \label{soln'}
da_{t}=-a_t\varsigma _{t}(h_{2}^{\top}) k^{+}k(t,Y_{t})\,(\mathrm{d}\tilde{W}%
_{t}-\varsigma_t(h_2)dt).
\end{equation}
But the fact that $a_t\equiv0$ is the unique solution of that equation is a consequence of Theorem 6 in Chapter V of \cite{protter2005stochastic}.
\end{proof}

We can now establish the final result of this section.
\begin{theorem}
\label{equivuniq} Let $\pi_0=\varsigma_0$ be a probability measure. Then
uniqueness of a measure valued solution   of  the evolution equation (\ref%
{eq:filterEqns:zakaifunctint}) is equivalent to uniqueness of a measure
valued solution   of the evolution equation (\ref{eq:filterEqns:ks}).
\end{theorem}

\begin{proof}
Let us assume first that there exists a unique solution of the evolution
equation (\ref{eq:filterEqns:ks}). Let $\pi ^{1}$ and $\pi ^{2}$ be two
solutions of the evolution equation (\ref{eq:filterEqns:zakaifunctint}) and
let $\pi ^{1}\left( 1\right) $ and $\pi ^{2}\left(1\right) $ be the
corresponding total mass processes. From (\ref{eq:filterEqns:zakaifunctint}%
), we deduce that these processes satisfy the evolution equation (take $%
\varphi _{t}\equiv 1$ in (\ref{eq:filterEqns:zakaifunctint})) 
\begin{equation}
\pi _{t}^{i}(1)=\pi _{0}^{i}(1)+\int_{0}^{t}\pi _{s}^{i}(h_{2}^{\top
})k^{+}\left( s,Y_{s}\right) \left( \mathrm{d}Y_{s}-h_{1}\left(
s,Y_{s}\right) \mathrm{d}s\right) .  \label{totalmass}
\end{equation}%
Define the normalized version of $\pi ^{1}$ and $\pi ^{2}$, $\varsigma
^{i}\equiv \frac{\pi ^{i}}{\pi ^{i}(1)}$. Then both $\varsigma ^{1}$ and $%
\varsigma ^{2}$ satisfy the evolution equation (\ref{eq:filterEqns:ks}), as
follows from the argument in the last part of the proof of Theorem \ref%
{th:ex}. It follows that $\varsigma ^{1}=\varsigma ^{2}=\varsigma $%
.~Moreover from equation (\ref{totalmass}), we deduce that 
\begin{eqnarray*}
\pi _{t}^{i}(1) &=&1+\int_{0}^{t}\pi _{s}^{i}(1)\frac{\pi
_{s}^{i}(h_{2}^{\top })}{\pi _{s}^{i}(1)}k^{+}\left( s,Y_{s}\right) dN_{t} \\
&=&1+\int_{0}^{t}\pi _{s}^{i}(1)\varsigma _{t}(h_{2}^{\top })k^{+}\left(
s,Y_{s}\right) dN_{t}
\end{eqnarray*}
In other words, both $\pi^{1}\left( 1\right) $ and $\pi^{2}\left( 1\right) $
are solutions of the equation (\ref{abstracttotalmass}) and therefore must
coincide by Lemma \ref{lemma:mass}. Hence 
\begin{equation*}
\pi ^{1}=\pi ^{1}\left( 1\right) \varsigma =\pi ^{2}\left( 1\right)
\varsigma =\pi ^{2}.
\end{equation*}%
Hence the evolution equation (\ref{eq:filterEqns:zakaifunctint}) has a
unique solution.

Now let us assume that the evolution equation (\ref%
{eq:filterEqns:zakaifunctint}) has a unique solution. Let $\varsigma ^{1}$
and $\varsigma ^{2}$ be two solutions of the evolution equation (\ref%
{eq:filterEqns:ks}) and consider $j^{\varsigma^{i}},$ $i=1,2$ the
corresponding solutions of the equation (\ref{abstracttotalmass}). Then $%
j^{\varsigma^1}\varsigma^{1}$ and $j^{\varsigma^{2}} \varsigma^{2}$ are
solutions of the evolution equation (\ref{eq:filterEqns:zakaifunctint}) by
Lemma \ref{lemma:prob} and therefore must be equal, $j^{\varsigma^1}%
\varsigma^{1}=j^{\varsigma^{2}} \varsigma^{2}=\pi $. It follows that%
\begin{equation*}
j^{\varsigma^1}=j^{\varsigma^1}\varsigma^{1}(1)=\pi \left( 1\right)
=j^{\varsigma^{2}}\varsigma^{2}\left( 1\right) =j^{\varsigma^{2}},
\end{equation*}
since {\ $\varsigma^{i}$ are probability measures}, again, by Lemma \ref%
{lemma:prob}. In other words for both $i=1,2$, $j^{\varsigma^{i}}$ coincides
with the total mass process of $\pi $. This gives us 
\begin{equation*}
\varsigma^{1}=\frac{\pi }{j^{\varsigma^1}}=\frac{\pi }{j^{\varsigma^2}}%
=\varsigma^{2}
\end{equation*}
and hence the evolution equation (\ref{eq:filterEqns:ks}) has a unique
solution.
\end{proof}

\section{The Backward SPDE approach to uniqueness}

Following from Theorem \ref{th:ex}, the process $\pi _{t}$ satisfies
the evolution equation 
\begin{equation*}
	d\pi _{t}(\varphi _{t})=\pi _{t}(\partial _{t}\varphi _{t}+A\varphi
	_{t})dt+\pi _{t}(B^{j}\varphi _{t})d\tilde{W}_{t}^{j}.
\end{equation*}%
for any function $\varphi \in C_{b}^{1,2}([0,\infty )\times \mathbb{R}%
^{d}))$, where we use the convention of summation over repeated indices, and the
notation
\begin{equation*}
B_t^{j}\varphi_t :=\left( \nabla \varphi_t [\bar{g}k^+k](t,\cdot,Y_t)%
\right)_{j}+([h_{2}^Tk^+k](t,\cdot,Y_t))_{j}\varphi_t ,\ 1\leq j\leq m.
\end{equation*}
 In the classical filtering problem (where $Y_t$ does not appear in the
coefficients), one associates to the Zakai equation a proper adjoint
backward PDE, which allows to establish
uniqueness of a measure valued solution to the Zakai equation (this is Bensoussan's approach to showing uniqueness of the solution of the  Zakai equation, see \cite{bensoussan} for details).

However, in our situation, this approach is not feasible because the backward partial differential equation would involve the observation process $Y_t$
in its coefficients, resulting in a solution that is not adapted at each time $t$
to the past of that process. To address this issue, we employ an adjoint backward stochastic partial differential equation (SPDE) instead of an adjoint backward PDE. The solution to the SPDE remains adapted at each time $t$ to the past of the observation process. This approach will be developed in the next section. To facilitate this, we establish a new type of It\^o formula, which is essential for leveraging the duality between the Zakai equation and an adjoint BSPDE.

Let $\mathcal{E}_{T}$ be the space of progressively measurable (with respect
to the augmented filtration generated by $\tilde{W}$) processes $\{u_t,\ 0\le t\le T\}$ such that 
\begin{equation*}
u\in C([0,T]; C^2_b(\mathbb{R}^d))
\end{equation*}
and 
\begin{equation}
u _{t}=u_{0}+\int_{0}^{t}\Sigma _{s}ds+\int_{0}^{t}\Lambda _{s}^{j}d\tilde{W}
_{s}^{j},\ 0\le t\le T\, ,  \label{identity}
\end{equation}
where $\Sigma ,\Lambda ^{j}$ are progressively measurable $C_{b}(\mathbb{R}%
^{d})$-valued processes such that 
\begin{equation*}
\Sigma\in L^1(0,T;C_b(\mathbb{R}^d)),\quad \Lambda^j\in L^2(0,T;C^2_b(%
\mathbb{R}^d))\, .
\end{equation*}

We denote by $\mathcal{U}$ the set of progressively measurable processes
with values in $\mathcal{M}_F(\mathbb{R}^d)$ which satisfy for all $T>0$ 
\begin{equation*}
\sup_{0\le t\le T}\mu_t(\mathbf{{1})<\infty \ \text{ a.s.}}
\end{equation*}
We shall say that $\mu\in\mathcal{U}$ solves the Zakai equation if for any $%
\varphi\in C^2_b(\mathbb{R}^d)$, 
\begin{equation*}
\mu_{t}(\varphi )=\mu_0(\varphi)+\int_0^t\mu_{s}(A_s\varphi
)ds+\int_0^t\mu_{s}(B^{j}_s\varphi )d\tilde{W}_{s}^{j},\ t\ge0\,.
\end{equation*}
We will now establish a useful It\^o formula.

\begin{theorem}
\label{Zakaiextended} For any $u \in \mathcal{E}_{T}$ of the form (\ref%
{identity}) and any $\mu \in \mathcal{U}$ that solves the Zakai equation we
have (again with the convention of summation over repeated indices) 
\begin{align}  \label{ito}
\mu _{t}\left( u _{t}\right) =\mu _{0}\left( u _{0}\right) +\int_{0}^{t}\mu
_{s}\left( A_{s}u _{s}+\Sigma _{s}+B_{s}^{j}\Lambda _{s}^{j}\right) ds
+\int_{0}^{t}\mu _{s}\left( B_{s}^{j}u _{s}+\Lambda _{s}^{j}\right) d\tilde{W%
}_{s}^{j}.
\end{align}
\end{theorem}

\begin{proof}
We fix $t>0$. For any $n\ge1$, $0\le s\le t$ and $x\in\mathbb{R}^d$, we let 
\begin{equation}  \label{approx}
\Lambda^j_n(s,x)=\sum_{p=0}^{n-1}\frac{n}{t}\int_{(\frac{p-1}{n}t)^+}^{\frac{%
p}{n}t} \Lambda^j(r,x)dr\mathbf{1}_{[\frac{p}{n}t,\frac{p+1}{n}t)}(s),
\end{equation}
and define 
\begin{equation*}
u_n(t,x)=u_0(x)+\int_0^t\Sigma_s(x)ds+\int_0^t\Lambda^j_n(s,x)d\tilde{W}%
_{s}^{j}\,.
\end{equation*}
It is easy to see that $\Lambda^j_n\to\Lambda^j$ in $L^2((0,T);C^2_b(\mathbb{%
R}^d))$ a.s.. Consequently $u_n\to u$ in $C([0,T];C^2_b(\mathbb{R}^d))$ in
probability. Hence if we show \eqref{ito} with $(u,\Lambda^j)$ replaced by $%
(u_n,\Lambda^j_n)$, the result will follow by taking the limit as $%
n\to\infty $. So from now on, we assume that $\Lambda^j=\Lambda^j_n$ is
given by \eqref{approx}, and delete the index $n$. Now it suffices to prove %
\eqref{ito} with $(0,t)$ replaced by $(a,b)$, with for some $1\le k\le n-1$, 
$\frac{k-1}{n}t\le a<b\le \frac{k}{n}t$. In other words, all we need to show
is that

\begin{align*}
\mu _{b}\left( u(b)\right) &=\mu _{a}\left( u (a)\right) +\int_{a}^{b}\mu
_{s}\left( A_{s}u(s)+\Sigma(s)+B_{s}^{j}\Lambda^{j}(a)\right) ds \\
&\quad+\int_{a}^{b}\mu _{s}\left( B_{s}^{j}u(s)+\Lambda^{j}(a)\right) d%
\tilde{W}_{s}^{j}.
\end{align*}
Let now $a=s_0<s_1<\cdots<s_{n^{\prime }}=b$, with $s_i=a+i\frac{b-a}{%
n^{\prime }}$, where $n^{\prime }$ is an integer which will eventually tend
to $+\infty$ (while $n$ is kept fixed). We have 
\begin{align*}
&\mu_{s_{i+1}} (u(s_{i+1}))-\mu_{s_i}(u(s_i)) \\
&= \mu_{s_{i+1}}(u(s_i))-\mu_{s_i}(u(s_i)) +\mu_{s_{i+1}}(u(s_{i+1}) -
u(s_i)) \\
&=\int_{s_i}^{s_{i+1}} \mu_s(A_s u(s_i))ds+\int_{s_i}^{s_{i+1}}
\mu_s(B^j_su(s_i))d\tilde{W}_{s}^{j} \\
&\quad+\int_{s_i}^{s_{i+1}}\mu_{s_{i+1}}(\Sigma(s))ds+\mu_{s_{i+1}}(%
\Lambda^j(a)) (\tilde{W}_{s_{i+1}}^{j}-\tilde{W}_{s_i}^{j}) \\
&=\int_{s_i}^{s_{i+1}} \mu_s(A_su(s_i))ds+\int_{s_i}^{s_{i+1}}
\mu_s(B^j_su(s_i))d\tilde{W}_{s}^{j} \\
&\quad+\int_{s_i}^{s_{i+1}}\mu_{s_{i+1}}(\Sigma(s))ds
+\mu_{s_{i}}(\Lambda^j(a)) (\tilde{W}_{s_{i+1}}^{j}-\tilde{W}_{s_i}^{j}) \\
&\quad+ (\tilde{W}_{s_{i+1}}^{j}-\tilde{W}_{s_i}^{j})\int_{s_i}^{s_{i+1}}%
\mu_s(A_s\Lambda^j(a))ds + (\tilde{W}_{s_{i+1}}^{j}-\tilde{W}%
_{s_i}^{j})\int_{s_i}^{s_{i+1}} \mu_s(B^{j^{\prime }}\Lambda^{j^{\prime
}}(a))d\tilde{W}^{j^{\prime }}_s,
\end{align*}
We wish to show that, as $n^{\prime }\to\infty$, 
\begin{align*}
\sum_{i=0}^{n^{\prime }-1}&\left[\mu_{s_{i+1}} (u(s_{i+1}))-\mu_{s_i}(u(s_i))%
\right] \\
&\to \int_a^b \mu_s(A_su(s)+\Sigma(s)+B^j_s\Lambda^j(a))
ds+\int_a^b\mu_s(B^j_su(s)+\Lambda^j(a))d\tilde W^j_s
\end{align*}
in probability.

We first show that in probability, as $n^{\prime }\to\infty$, 
\begin{equation*}
\sum_{i=0}^{n^{\prime }-1}\left[\int_{s_i}^{s_{i+1}}
\mu_s(A_su(s_i))ds+\int_{s_i}^{s_{i+1}} \mu_s(B^j_su(s_i))d\tilde{W}_{s}^{j} %
\right] \to \int_a^b\mu_s(A_su(s))ds+\int_a^b \mu_s(B^j_su(s))d\tilde{W}_s^j
\end{equation*}
This statement follows from the fact that, since $u\in C([0,T];C^2_b(\mathbb{%
R}^d))$, and 
\begin{align*}
\sup_i\sup_{s_i\le s<s_{i+1}}|\mu_s(A_s [u(s_i)-u(s)])&\le \sup_s\mu_s(%
\mathbf{1})\sup_i\sup_{s_i\le s<s_{i+1}}\|u(s_i)-u(s)\|_{C^2_b(\mathbb{R}%
^r)}, \\
\sup_i\sup_{s_i\le s<s_{i+1}}|\mu_s(B^j_s [u(s_i)-u(s)])&\le \sup_s\mu_s(%
\mathbf{1})\sup_i\sup_{s_i\le s<s_{i+1}}\|u(s_i)-u(s)\|_{C^1_b(\mathbb{R}%
^r)},\ 1\le j\le \ell^{\prime }\,.
\end{align*}
we have that in probability, as $n^{\prime }\to\infty$, 
\begin{align*}
\sum_{i=0}^{n^{\prime }-1}\mathbf{1}_{[s_i,s_{i+1})}\mu_s(A_su(s_i))&\to
\mu_s(A_su(s)) \quad\text{in }C([0,T]), \\
\sum_{i=0}^{n^{\prime }-1}\mathbf{1}_{[s_i,s_{i+1})}\mu_s(B^j_su(s_i))&\to
\mu_s(B^j_su(s)) \quad\text{in }C([0,T]),\ 1\le j\le \ell^{\prime }\,.
\end{align*}

Secondly, since $\Lambda^j(a)\in C^2_b(\mathbb{R}^d))$, it follows from the
Zakai equation that $s\to\mu_s(\Lambda^j(a))$ is continuous. 
Hence clearly 
\begin{equation*}
\sum_{i=0}^{n^{\prime }-1}\mu_{s_{i}}(\Lambda^j(a)) (\tilde{W}_{s_{i+1}}^{j}-%
\tilde{W}_{s_i}^{j}) \to\int_a^b\mu_s(\Lambda^j(a)) d\tilde{W}_s\,.
\end{equation*}
Moreover, we know that $\mu_\cdot(B_j\Lambda^j(a))\in L^\infty(0,T)$. Hence
a classical arguments yields that 
\begin{equation*}
\sum_{i=0}^{n^{\prime }-1}(\tilde{W}_{s_{i+1}}^{j}-\tilde{W}%
_{s_i}^{j})\int_{s_i}^{s_{i+1}} \mu_s(B^k\Lambda^k(a))d\tilde{W}%
^k_s\to\int_a^b\mu_s(B^j\Lambda^j(a))ds\,.
\end{equation*}
Indeed, the limit is the joint quadratic variation on the interval $[a,b]$
of the two martingales $\tilde{W}_t$ and $\int_0^t\mu_s(B^k\Lambda^k(a))d%
\tilde{W}^k_s$.

We also note that 
\begin{align*}
\left|\sum_{i=0}^{n^{\prime }-1}(\tilde{W}_{s_{i+1}}^{j}-\tilde{W}%
_{s_i}^{j})\int_{s_i}^{s_{i+1}}\mu_s(A_s\Lambda^j(a))ds\right| &\le\sup_i|%
\tilde{W}_{s_{i+1}}^{j}-\tilde{W}_{s_i}^{j}|\int_a^b|\mu_s(A_s%
\Lambda^j(a))|ds \\
&\le C\sup_i|\tilde{W}_{s_{i+1}}^{j}-\tilde{W}_{s_i}^{j}|\int_a^b|\mu_s(%
\mathbf{1})|ds \\
&\to0\,.
\end{align*}
Finally, we show that 
\begin{equation}  \label{sigma}
\sum_{i=0}^{n^{\prime
}-1}\int_{s_i}^{s_{i+1}}\mu_{s_{i+1}}(\Sigma(s))ds\to\int_a^b\mu_s(%
\Sigma(s))ds\,.
\end{equation}
We approximate $\Sigma$ in $L^1(0,T;C_b(\mathbb{R}^d))$ by a sequence in $%
C([0,T];C^2_b(\mathbb{R}^d))$. For each $\varepsilon>0$, let $%
\Sigma_\varepsilon\in C([0,T];C^2_b(\mathbb{R}^d))$ be such that $%
\int_0^T\sup_x|\Sigma(t,x)-\Sigma_\varepsilon(t,x)|dt\le \varepsilon$. We
have 
\begin{align*}
\left|\sum_{i=0}^{n^{\prime
}-1}\int_{s_i}^{s_{i+1}}\mu_{s_{i+1}}(\Sigma(s))ds-\int_a^b\mu_s(%
\Sigma(s))ds\right|&\le \left|\sum_{i=0}^{n^{\prime
}-1}\int_{s_i}^{s_{i+1}}\mu_{s_{i+1}}(\Sigma(s))ds-\sum_{i=0}^{n^{\prime
}-1}\int_{s_i}^{s_{i+1}}\mu_{s_{i+1}}(\Sigma_\varepsilon(s))ds\right| \\
&\quad +\left|\sum_{i=0}^{n^{\prime
}-1}\int_{s_i}^{s_{i+1}}[\mu_{s_{i+1}}(\Sigma_\varepsilon(s))-\mu_s(\Sigma_%
\varepsilon(s))]ds\right| \\
&\quad+\left|\int_a^b\mu_s(\Sigma_\varepsilon(s))ds-\int_a^b\mu_s(%
\Sigma(s))ds\right|
\end{align*}
We observe that the first and the last term on the right hand side of the above inequality are bounded by $(b-a)\varepsilon%
\sup_{a\le s\le b}\mu_s(\mathbf{1})$. It thus remains to show that for each $%
\varepsilon>0$ fixed, the second term tends to $0$, as $n^{\prime }\to\infty$%
. This follows from the fact that, since $\mu_t$ solves the Zakai equation,
for any $s_i\le s\le s_{i+1}$, 
\begin{align*}
\mu_{s_{i+1}}(\Sigma_\varepsilon(s))-\mu_{s}(\Sigma_\varepsilon(s))=%
\int_{s}^{s_{i+1}}\mu_r(A_r\Sigma_\varepsilon(s))dr
+\int_{s}^{s_{i+1}}\mu_r(B^j_r\Sigma_\varepsilon(s))d\tilde{W}^j_r\, .
\end{align*}
We first note that 
\begin{align*}
\left|\sum_{i=0}^{n^{\prime
}-1}\int_{s_i}^{s_{i+1}}\int_{s}^{s_{i+1}}\mu_r(A_r\Sigma_%
\varepsilon(s))drds\right| \le\sup_{a\le s\le b}\left(\mu_s(\mathbf{1}%
)\|\Sigma_\varepsilon(s)\|_{C^2_b}\right)\times\sum_i(s_{i+1}-s_i)^2/2,
\end{align*}
which tends to $0$, as $n^{\prime }\to\infty$, for any $\varepsilon>0$
fixed. Moreover, for any $M>0$, $\delta>0$, 
\begin{align*}
\P \left(\left|\sum_{i=0}^{n^{\prime
}-1}\int_{s_i}^{s_{i+1}}ds\int_{s}^{s_{i+1}}\mu_r(B^j_r\Sigma_%
\varepsilon(s))d\tilde{W}^j_r\right|>\delta\right)&\le \P \left(\sup_{a\le
s\le b}(\mu_s(\mathbf{1})\|\Sigma_\varepsilon(s)\|_{C^1_b})>M\right) \\
&\quad+Ck\frac{M}{\delta}\sum_i(s_{i+1}-s_i)^{3/2}
\end{align*}
For any $M>0$ and $\delta>0$ fixed, the second term on the right tends to $0$
as $n^{\prime }\to\infty$, while the first term tends to $0$ as $M\to\infty$%
, with $\varepsilon>0$ fixed. \eqref{sigma} has been established.
\end{proof}

We will establish the above result with the same assumptions on $\Sigma$,
and $\Lambda^j\in L^2((0,T);C^1_b(\mathbb{R}^d))$. However, the processes $u$%
, $\Sigma$ and $\Lambda^j$ will be given a Sobolev--space valued processus,
and the fact that they take their values in $C^2_b(\mathbb{R}^d)$, $C_b(%
\mathbb{R}^d)$ and $C^1_b(\mathbb{R}^d)$ respectively, will be a consequence
of classical Sobolev embedding theorems. Hence the result which will be
useful to us is the following theorem, where $H^m:=H^m(\mathbb{R}^d)$ (with $%
m\ge0$ an integer) denotes the Sobolev space of square integrable functions
whose distributional derivatives up to order $m$ are all square integrable.
In particular, $H^0=L^2(\mathbb{R}^d)$. Moreover we will denote by $\mathcal{%
P}$ the $\sigma$-field of predictable subsets of $\mathbb{R}_+\times\Omega$,
and for any Hilbert space $H$, $L^2_\mathcal{P}(\Omega\times[0,T];H)$ will
denote the set of $H$ valued processes which are $\mathcal{P}$ measurable,
and square integrable with respect to the product measure $dt\times d\P $.


For the proof of Theorem \ref{Itoform} below, we need the following
fundamental result on SPDEs:

\begin{proposition}
\label{SPDE-EP} Let for some $m\ge0$ $u_0\in L^2(\Omega;H^m)$, $f\in L^2_%
\mathcal{P}(\Omega\times[0,T];H^{m-1})$ and for $1\le j\le \ell^{\prime }$, $%
g^j\in L^2_\mathcal{P}(\Omega\times[0,T];H^{m})$. Then the SPDE 
\begin{equation*}
 u_t=u_0+\int_0^t [\Delta u_s+f_s]ds+\int_0^t g^j_sdW^j_s,\ t\ge0
\end{equation*}
has a unique solution $u\in L^2_\mathcal{P}(\Omega\times[0,T];H^{m+1})\cap
L^2(\Omega;C([0,T];H^m))$. Moreover, the mapping $(f,g^1,\ldots,g^k)\mapsto
u $ is continuous from $L^2_\mathcal{P}(\Omega\times[0,T];H^{m-1})\times%
\left(L^2_\mathcal{P}(\Omega\times[0,T];H^{m}\right)^{\ell^{\prime }}$ into $%
L^2_\mathcal{P}(\Omega\times[0,T];H^{m+1})\cap L^2(\Omega;C([0,T];H^m))$.
\end{proposition}

\begin{proof}
The existence and uniqueness result in the case $m=0$ is a particular case
of Theorem 1.4 in \cite{Pard-Stoch} (see also \cite{pp}). The continuity
follows readily from the estimates there. The result in the case $m\ge1$ is
deduced as follows. For any $1\le i\le d$, $v_i:=\partial u/\partial x_i$ is
the solution of an equation to which the result for $m=0$ can be applied.
This establishes the result for $m=1$. The result for $m>1$ is
obtained inductively by taking higher order derivatives.
\end{proof}

\begin{theorem}
\label{Itoform} Suppose that, for some $m>d/2+2$, $u\in L^2_\mathcal{P}
(\Omega\times[0,T];H^{m}))$ and moreover for any $0\le t\le T$, 
\begin{equation*}
u(t)=u_0+\int_0^t\Sigma(s)ds+\int_0^t\Lambda^j(s)d\tilde{W}^j(s),
\end{equation*}
where $u_0\in L^2(\Omega;H^m)$ is $\mathcal{F}_0$ measurable, $\Sigma\in L^2_%
\mathcal{P}(\Omega\times(0,T),H^{m-2})$ and $\Lambda^j\in L^2_\mathcal{P}%
(\Omega\times(0,T);H^{m-1})$ for all $1\le j\le \ell^{\prime }$. Then the
It\^o formula in Theorem \ref{Zakaiextended} still holds, i.e. 
\begin{align*}
\mu _{t}\left( u _{t}\right) =\mu _{0}\left( u _{0}\right) +\int_{0}^{t}\mu
_{s}\left( A_{s}u _{s}+\Sigma _{s}+B_{s}^{j}\Lambda _{s}^{j}\right) ds
+\int_{0}^{t}\mu _{s}\left( B_{s}^{j}u _{s}+\Lambda _{s}^{j}\right) d\tilde{W%
}_{s}^{j}.
\end{align*}
\end{theorem}

\begin{proof}
Let for all $1\le j\le \ell^{\prime }$ $\{\Lambda^j_n,\ n\ge1\}$ denote a
sequence in $L^2_\mathcal{P}(\Omega\times(0,T);H^{m})$, such that $%
\Lambda^j_n\to\Lambda^j$ in $L^2_\mathcal{P}(\Omega\times(0,T);H^{m-1})$.
Let moreover $\{u_n,\ n\ge1\}$ (resp. $\{\Sigma_n,\ n\ge1\}$) denote a
sequence in $L^2_\mathcal{P}(\Omega\times(0,T),H^{m+1})$ (resp. in $L^2_%
\mathcal{P}(\Omega\times(0,T),H^{m-1})$), such that $u_n\to u$ in $L^2_%
\mathcal{P}(\Omega\times(0,T),H^{m})$ (resp. $\Sigma_n\to\Sigma$ in $L^2_%
\mathcal{P}(\Omega\times(0,T),H^{m-2})$). We now define for each $n\ge1$ $%
v_n $ as the solution of the following SPDE (where $\Delta$ denotes the
Laplace operator) : 
\begin{equation*}
v_n(t)=u_0+\int_0^t[\Delta v_n(s)-\Delta
u_n(s)+\Sigma_n(s)]ds+\int_0^t\Lambda^j_n(s)d\tilde{W}^j(s),\ 0\le t\le T\, .
\end{equation*}
We shall now use repeatedly the results in Proposition \ref{SPDE-EP}. It is
plain that this SPDE has a unique solution $v_n\in L^2_\mathcal{P}(\Omega\times(0,T);H^{m+1})\cap L^2_\mathcal{P}(\Omega;C([0,T];H^m))$, 
and
as $n\to\infty$, $v_n\to\tilde{u}$ in $L^2_\mathcal{P}(\Omega\times(0,T);H^{m})$, 
where $\tilde{u}$ is the unique solution in $L^2(\Omega\times(0,T);H^{m})$ of the SPDE 
\begin{equation*}
\tilde{u}(t)=u_0+\int_0^t[\Delta \tilde{u}(s)-\Delta
u(s)+\Sigma(s)]ds+\int_0^t\Lambda^j(s)d\tilde{W}^j(s),\ 0\le t\le T\, .
\end{equation*}
But $u\in L^2(\Omega\times(0,T);H^{m})$ is a solution of that equation.
Hence $\tilde{u}=u$. Now from Theorem \ref{Zakaiextended}, which we can use
thanks to the Sobolev embedding, which in particular tells us that $H^m\subset C^2_b(\mathbb{R}^d)$, 
\begin{align*}
\mu _{t}\left( v _n({t})\right) =&\mu _{0}\left( u _{0}\right)
+\int_{0}^{t}\mu _{s}\left( A_{s}v _n(s)+\Sigma _n({s})+
\Delta(v_n-u_n)(s)+B_{s}^{j}\Lambda _{n}^{j}(s)\right) ds \\
&+\int_{0}^{t}\mu _{s}\left( B_{s}^{j}v_n(s)+\Lambda^{j} _n(s)\right) d\tilde{W}_{s}^{j}.
\end{align*}
Now we can take the limit in that identity as $n\to\infty$, which
yields the result. Indeed, as $n\to\infty$, for any $t>0$, $v_n(t)\to u(t)$ in $L^2(\Omega;H^{m-1})$, 
\begin{equation*}
Av_n+\Sigma_n+\Delta(v_n-u_n)+B^j\Lambda^j_n\to Au+\Sigma+B^j\Lambda^j\ 
\text{ in }L^2_\mathcal{P}(\Omega\times(0,T);H^{m-2})
\end{equation*}
 for $1\le j\le \ell^{\prime }$, 
\begin{equation*}
B^jv_n+\Lambda^j_n\to B^ju+\Lambda^j\ \text{ in }L^2_\mathcal{P}%
(\Omega\times(0,T);H^{m-1}),
\end{equation*}
$H^{m-2}(\mathbb{R}^d)\subset C_b(\mathbb{R}^d)$ with continuous injection,
and $\sup_{0\le t\le T}\mu_t(\mathbf{1})<\infty$. Combining those facts, we deduce that that as $n\to\infty$, the following convergences hold in probability: 
\begin{align*}
\mu_t(v_n(t))&\to \mu_t(u(t)), \\
\mu_s(A_sv_n(s)+\Sigma_s+\Delta(v_n-u_n)(s)+B^j_s\Lambda^j_n(s))&\to
\mu_s(A_su_s+\Sigma_s+B^j_s\Lambda^j_s)\ \text{ in }L^2(0,T)
\end{align*}
and for $1\le j\le k$, 
\begin{equation*}
\mu_s(B^j_sv_n(s)+\Lambda^j_n(s))\to \mu_s(B^j_su_s+\Lambda^j_s)\ \text{ in }%
L^2(0,T)\, .
\end{equation*}
The result follows.
\end{proof}

\section{A system of BSPDEs}

In the following we will make use of a complex valued $u \in \mathcal{E}_{T}$
of the form (\ref{identity}), which will be the solution of the BSPDE 
\begin{equation}  \label{complexBSPDE}
du _{t}=-\left( Au _{t}+B^{j}v_{t}^{j}+ir_{t}^{j}B^{j}u
_{t}+ir_{t}^{j}v_{t}^{j}\right) dt+v_{t}^{j}dW_{t}^{j}\,,\quad u
_{T}=\varphi,
\end{equation}
where again we adopt the convention of summation of the repeated index $j$
from $j=1$ to $j=\ell^{\prime }$. Under the Assumptions \textbf{AA}$_m$ to be
specified below, as a result of Theorem \ref{bspde}, $u$ will satisfy the
assumptions of Theorem \ref{Itoform}.

We write below the corresponding equations of the real, respectively, the
imaginary part of $u$. In other words, assume that $u =u ^{1}+iu ^{2}$. Then 
$\left( u ^{1},u ^{2}\right) $ satisfy the following system of BSPDEs%
\begin{equation}  \label{SBSPDE}
\begin{split}
du _{t}^{1} &=-\left(
Au_{t}^{1}+B^{j}v_{t}^{1,j}-r_{t}^{j}B^{j}u_{t}^{2}-r_{t}^{j}v_{t}^{2,j}%
\right) dt+v_{t}^{1,j}dW_{t}^{j}\,,\quad u _{T}^{1}=\varphi, \\
du _{t}^{2} &=-\left(
Au_{t}^{2}+B^{j}v_{t}^{2,j}+r_{t}^{j}B^{j}u_{t}^{1}+r_{t}^{j}v_{t}^{1,j}%
\right) dt+v_{t}^{2,j}dW_{t}^{j}\,,\quad u _{T}^{2}=0.
\end{split}%
\end{equation}
We need to extend to the above system of BSPDEs the results from Du and Meng 
\cite{DM} and from Du, Tang and Zhang \cite{DTZ}, which are established for
a single BSPDE, of the same type. Note that the factor of $dt$ in the $u^1$
(resp. $u^2$) equation involves second order derivatives of $u^1$ (resp. $%
u^2 $) and first order derivatives of $u^2$ (resp. $u^1$), together with
first order derivatives of $v^1$ (resp. $v^2$) and zero-th order derivatives
of $v^2$ (resp. $v^1$). Hence the coupling between the two BSPDEs comes
through terms of lower order, which is essential for our extension from the
results for a single BSPDE to work.

We now first state and prove the extension of Theorem 2.3 from \cite{DM} to
our system. 

\begin{theorem}
\label{th:SP} \label{firsttheorem} In addition to the assumptions \textbf{E}
and \textbf{U}, let us suppose that for some $\kappa>0$, 
\begin{equation}  \label{coercivity}
gg^{T}(t,x)\ge\kappa I,\ \forall (t,x)\in[0,T]\times\mathbb{R}^d,\ \text{a.s.%
},
\end{equation}
and for some integer $n\ge1$, any multi-index $\alpha$ with $|\alpha|\le n$, 
\begin{equation}  \label{smooth-coeff}
\text{ess sup}_{\Omega\times[0,T]\times\mathbb{R}^d}(|D^\alpha a|+|D^\alpha
f|+|D^\alpha\bar{g}|+|D^\alpha h_2|)\le K.
\end{equation}
Finally we assume that $\varphi\in H^{n+1}$.

Then the system of BSPDEs \eqref{SBSPDE} has a unique solution such that for 
$i=1,2$, 
\begin{equation*}
u^i\in L^2_\mathcal{P}(\Omega\times[0,T];H^{n+2})\cap
L^2(\Omega;L^\infty([0,T];H^{n+1})), v^i\in L^2_\mathcal{P}(\Omega\times[0,T]%
;(H^{n+1})^{\otimes\ell^{\prime }}).
\end{equation*}
\end{theorem}

\begin{proof}
We first need to extend Proposition 3.2 together with Theorem 2.1 from \cite%
{DM}. Let $V:=H^1(\mathbb{R}^d)\times H^1(\mathbb{R}^d)$, $H=L^2(\mathbb{R}%
^d)\times L^2(\mathbb{R}^d)$, so that if we identify $H$ with its dual, $%
V^{\prime }$ is identified with $H^{-1}(\mathbb{R}^d)\times H^{-1}(\mathbb{R}%
^d)$. Referring to the notations in the proof of Proposition 3.2 in \cite{DM}%
, we let 
\begin{equation*}
\mathcal{L}=%
\begin{pmatrix}
A & -r^jB^j \\ 
r^jB^j & A%
\end{pmatrix}
,\quad \mathcal{M}^j= 
\begin{pmatrix}
B^j & -r^j \\ 
r^j & B^j%
\end{pmatrix}
\,.
\end{equation*}
It is not hard to deduce from condition \eqref{coercivity} that Assumption
3.1 in \cite{DM} is satisfied, namely there exists $\lambda, C>0$ such that
for any $u\in V$, 
\begin{equation*}
2\langle u,\mathcal{L}u\rangle +\sum_{j=1}^{\ell^{\prime }} \|(\mathcal{M}%
^j)^Tu\|^2_{H}\le -\lambda\|u\|^2_V+C\| u\|^2_H,\ \|\mathcal{L}%
u\|_{V^{\prime }}\le C\|u\|_V\,.
\end{equation*}
It follows that the proofs of Proposition 3.2 and Theorem 2.1 from \cite{DM}
are easily adapted to yield that provided the assumption \eqref{smooth-coeff}
is satisfied with $n=0$ and $\varphi\in L^2(\mathbb{R}^d)$, our system has a
unique solution such that $(u^i,v^i)\in L^2(\Omega\times(0,T);H^1(\mathbb{R}%
^d)\times (L^2(\mathbb{R}^d))^{\otimes\ell^{\prime }}) $, $i=1,2$.

The rest of the proof follows exactly the lines of arguments in \cite{DM},
with obvious adaptations.
\end{proof}

%

Assumptions \textbf{AA}$_m$:

\begin{itemize}
\item \emph{(smoothness and boundedness of the coefficients)} All
coefficients are functions of $(\omega ,t,x)\in \Omega \times \left[ 0,T%
\right] \times \mathbb{R}^{d}$, which are $\mathcal{P}\otimes\mathcal{B}_d$
measurable and, for some integer $m$, (\ref{smooth-coeff}) is satisfied for
any multi-index $\alpha $ with $|\alpha |\leq \sup\{1,m\}$, and for $%
|\alpha|\le\sup\{2,m\}$ concerning the coefficients of the matrix $a$.


\item \emph{(Smoothness of the final condition)} The function $\varphi\in
H^{m}$.
\end{itemize}

In the following we will use the notation 
\begin{align*}
Q_{t}^{1,u,v}
&=Au_{t}^{1}+B^{j}v_{t}^{1,j}-r_{t}^{j}B^{j}u_{t}^{2}-r_{t}^{j}v_{t}^{2,j},
\\
Q_{t}^{2,u,v}
&=Au_{t}^{2}+B^{j}v_{t}^{2,j}+r_{t}^{j}B^{j}u_{t}^{1}+r_{t}^{j}v_{t}^{1,j}\,.
\end{align*}

\begin{theorem}
\label{bspde} Assume that for some integer $m$, \textbf{AA}$_m$ is
satisfied. Then the system of BSPDEs 
\begin{equation}  \label{splusf}
\begin{split}
du_{t}^{1} &=-Q_{t}^{1,u,v}dt+v_{t}^{1,j}dW_{t}^{j}\,, \\
du_{t}^{2} &=-Q_{t}^{2,u,v}dt+v_{t}^{2,j}dW_{t}^{j}\,.
\end{split}%
\end{equation}
with $(u_{T}^{1},u_{T}^{2})=(\varphi,0)$ has a solution $\left( (u^1,v^1),
(u^2,v^2)\right) $ such that for all $T>0$, $u^i\in L^2_\mathcal{P}%
(\Omega;C_w([0,T];H^{m}))$, $v^i\in L^2_\mathcal{P}(\Omega\times[0,T]%
;(H^{m-1})^{\otimes \ell^{\prime }})$, $i=1,2$ and we have for $i=1,2$ (here 
$\|\cdot\|_m$ stands for the norm in $H^m$ and $|\|\cdot|\|_m$ for the norm
in $(H^m)^{\otimes\ell^{\prime }}$) 
\begin{align}
E\left[ \sup_{t\leq T}\left\vert \left\vert u^i_{t}\right\vert
\right\vert_{m}^{2} \right]+ \mathbb{E}\left[ \int_{0}^{T}\left\vert
\left\vert \left\vert v^i_{t}+k^+k\bar{g}^T\nabla u^i_t\right\vert
\right\vert\right\vert _{m}^{2}\right]\leq &CE\left\vert \left\vert
\varphi\right\vert \right\vert _{m}^{2}  \label{estim1} \\
E\left[ \sup_{t\leq T}\left\vert \left\vert u^i_{t}\right\vert \right\vert
_{m}^{2}\right]+ E\left[ \int_{0}^{T}\left\vert \left\vert \left\vert
v^i_{t}\right\vert \right\vert\right\vert _{m-1}^{2}\right]\leq
&CE\left\vert \left\vert \varphi\right\vert \right\vert _{m}^{2}.
\label{estim2}
\end{align}

Moreover, if $m>d/2$ then $u$ is jointly continuous in $\left( t,x\right) $
almost surely and if $m >d/2+2$ then $\left( u,v\right)$ is a classical
solution of (\ref{splusf}). 
\end{theorem}

\begin{proof}
The proof is almost identical to the proof of Theorem 2.1 and Corollary 2.2
in \cite{DTZ}. Let us first explain how \eqref{estim2} follows from %
\eqref{estim1}. We first deduce from \eqref{estim1} that 
\begin{equation*}
E\left[ \sup_{t\leq T}\left\vert \left\vert u^i_{t}\right\vert
\right\vert_{m}^{2}\right]\leq CE\left\vert \left\vert \varphi\right\vert
\right\vert _{m}^{2}\,.
\end{equation*}
This clearly implies that 
\begin{equation*}
\mathbb{E}\int_0^T\left\vert \left\vert \left\vert k^+k\bar{g}^T\nabla
u^i_t\right\vert \right\vert\right\vert_{m-1}^{2}\leq CE\left\vert
\left\vert \varphi\right\vert \right\vert _{m}^{2}\,.
\end{equation*}
But \eqref{estim1} is also true for $m$ replaced by $m-1$, which implies
that 
\begin{equation*}
E \int_{0}^{T}\int_{0}^{T}\left\vert \left\vert \left\vert v^i_{t}+k^+k\bar{g%
}^T\nabla u^i_t\right\vert \right\vert\right\vert _{m-1}^{2}\leq
CE\left\vert \left\vert \varphi\right\vert \right\vert _{m-1}^{2},
\end{equation*}
and \eqref{estim2} now follows from these three inequalities.

Let us now explain how we obtain \eqref{estim1} in the simple case where $%
m=0 $. Below $\|\cdot\|$ (resp. $|\|\cdot|\|$) stands for $\|\cdot\|_0$
(resp. $|\|\cdot|\|_0$), and $(\cdot,\cdot)$ stands for the scalar product
in $H^0=L^2(\mathbb{R}^d)$. Suppose we have a smooth enough solution of %
\eqref{splusf}. Applying It\^o's formula to compute $\|u^i_t\|^2$ and
summing up the results for $i=1$ and $2$, we obtain 
\begin{align}  \label{energyequ}
\|\varphi\|^2=\|u^1_t\|^2+\|u^2_t\|^2+\int_t^T\{\||v^1_s|\|^2+\||v^2_s|%
\|^2-2(Q^{i,u,v}_s,u^i_s)\}ds +2\int_t^T(u^i_s,v^{i,j}_s)dW^j_s,
\end{align}
Then 
\begin{eqnarray*}
-2(Q_{s}^{i,u,v},u_{s}^{i})+\left( \Vert |v_{s}^{1}|\Vert ^{2}+\Vert
|v_{s}^{2}|\Vert ^{2}\right) &=&-2\left(
Au_{s}^{1}+B^{j}v_{s}^{1,j}-r_{s}^{j}B^{j}u_{s}^{2}-r_{s}^{j}v_{s}^{2,j},u_{s}^{1}\right)
\\
&&-2\left(
Au_{s}^{2}+B^{j}v_{s}^{2,j}+r_{s}^{j}B^{j}u_{s}^{1}+r_{s}^{j}v_{s}^{1,j},u_{s}^{2}\right)
\\
&&+\left( \Vert |v_{s}^{1}|\Vert ^{2}+\Vert |v_{s}^{2}|\Vert ^{2}\right) \\
&=&-2\left( Au_{s}^{1},u_{s}^{1}\right) -2\left( Au_{s}^{2},u_{s}^{2}\right)
\\
&&-2\left( B^{j}v_{s}^{1,j},u_{s}^{1}\right) -2\left(
B^{j}v_{t}^{2,j},u_{s}^{2}\right) \\
&&+2r_{s}^{j}\left( B^{j}u_{s}^{2},u_{s}^{1}\right) -2r_{s}^{j}\left(
B^{j}u_{s}^{1},u_{s}^{2}\right) +2r_{s}^{j}\left(
v_{s}^{2,j},u_{s}^{1}\right) -2r_{s}^{j}\left( v_{s}^{1,j},u_{s}^{2}\right)
\\
&&+\left( \Vert |v_{s}^{1}|\Vert ^{2}+\Vert |v_{s}^{2}|\Vert ^{2}\right).
\end{eqnarray*}%
By integration by parts, we deduce that 
\begin{eqnarray*}
&&-2(Q_{s}^{i,u,v},u_{s}^{i})+\left( \Vert |v_{s}^{1}|\Vert ^{2}+\Vert
|v_{s}^{2}|\Vert ^{2}\right) \\
&&\hspace{3cm}= ([\text{div}f-\frac{1}{2}D^{2}a]u_s^{i},u_s^{i})+
(a^{j^{\prime }l}\partial^{j^{\prime }}u_s^i,\partial^{l}u_s^j ) \\
&&\hspace{3.5cm}-2\left( \left( \nabla v_{s}^{1,j}[\bar{g}k^{+}k](t,\cdot
,Y_{t})\right) _{j},u_{s}^{1}\right) -2\left( ([h_{2}^{T}k^{+}k](t,\cdot
,Y_{t}))_{j}v_{s}^{1,j},u_{s}^{1}\right) \\
&&\hspace{3.5cm}-2\left( \left( \nabla v_{s}^{2,j}[\bar{g}k^{+}k](t,\cdot
,Y_{t})\right) _{j},u_{s}^{2}\right) -2\left( ([h_{2}^{T}k^{+}k](t,\cdot
,Y_{t}))_{j}v_{s}^{2,j},u_{s}^{2}\right) \\
&&\hspace{3.5cm}+2r_{s}^{j}\left( \left( \nabla u_{s}^{2}[\bar{g}%
k^{+}k](t,\cdot ,Y_{t})\right) _{j},u_{s}^{1}\right) +2r_{s}^{j}\left(
([h_{2}^{T}k^{+}k](t,\cdot ,Y_{t}))_{j}u_{s}^{2},u_{s}^{1}\right) \\
&&\hspace{3.5cm}-2r_{s}^{j}\left( \left( \nabla u_{s}^{1}[\bar{g}%
k^{+}k](t,\cdot ,Y_{t})\right) _{j},u_{s}^{2}\right) -2r_{s}^{j}\left(
([h_{2}^{T}k^{+}k](t,\cdot ,Y_{t}))_{j}u_{s}^{1},u_{s}^{2}\right) \\
&&\hspace{3.5cm}+2r_{s}^{j}\left( v_{s}^{2,j},u_{s}^{1}\right)
-2r_{s}^{j}\left( v_{s}^{1,j},u_{s}^{2}\right) \\
&&\hspace{3.5cm}+\left( \Vert |v_{s}^{1}|\Vert ^{2}+\Vert |v_{s}^{2}|\Vert
^{2}\right) \\
&&\hspace{3cm}\geq ([\text{div}f-\frac{1}{2}D^{2}a]u^{i},u^{i})+(gg^{T}%
\nabla u_{s}^{i},\nabla u_{s}^{i})+\sum_{i=1}^{2}|\Vert v_{s}^{i}+k^{+}k\bar{%
g}^{T}\nabla u_{s}^{i}|\Vert ^{2} \\
&&\hspace{3cm}+(\alpha ^{j}v^{i,j},u^{i}) \\
&&\hspace{3cm}+2r_{s}^{j}\left( \left( \nabla u_{s}^{2}[\bar{g}%
k^{+}k](t,\cdot ,Y_{t})\right) _{j},u_{s}^{1}\right) -2r_{s}^{j}\left(
\left( \nabla u_{s}^{1}[\bar{g}k^{+}k](t,\cdot ,Y_{t})\right)
_{j},u_{s}^{2}\right) \\
&&\hspace{3cm}+2r_{s}^{j}\left( v_{s}^{2,j},u_{s}^{1}\right)
-2r_{s}^{j}\left( v_{s}^{1,j},u_{s}^{2}\right),
\end{eqnarray*}%
where 
\begin{equation*}
\alpha ^{j}=2\sum_{p}\left( \partial _{p}[\bar{g}k^{+}k](t,\cdot
,Y_{t})\right) _{j}-2([h_{2}^{T}k^{+}k](t,\cdot ,Y_{t}))_{j}.
\end{equation*}%
In the above we have used the convention of summation over repeated indices and the
notation $\text{div}f= f^i_{x_i}$, and $D^2a:=\frac{\partial^2}{\partial
x_i\partial x_j}a^{i,j}$ (here as an exception the repeated indices $i$ and $%
j$ are both summed from $1$ to $d$), and the fact that $k^+k$ is a
projection operator, hence $\bar{g}\bar{g}^T\ge\bar{g}k^+k\bar{g}^T$, that 
\begin{align}
-2(Q^{i,u,v}_s,u^i_s)+\||v^1_s|\|^2+\||v^2_s|\|^2&\ge(gg^T\nabla
u^i_s,\nabla u^i_s)+([\text{div} f-\frac{1}{2}D^2a]u^i,u^i)
+\sum_{i=1}^2|\|v^i_s+k^+k\bar{g}^T\nabla u^i_s|\|^2  \notag \\
&\quad +(\alpha^{i,j}v^{i,j},u^i)+(-1)^{3-i}2r^j_s((k^+k\bar{g}^T)^j\nabla
u^{3-i}_s+v^{3-i,j}_s,u^i_s)  \notag \\
&\ge\sum_{i=1}^2\left\{\frac{1}{2} |\|v^i_s+k^+k\bar{g}^T\nabla u^i_s|\|^2-C
\|u^i_s\|^2\right\} ,  \label{ineq}
\end{align}
Let us justify the last inequality. Since $r^j_s$ is bounded, there exists a
constant $C$ such that 
\begin{equation*}
(-1)^{3-i}2r^j_s((k^+k\bar{g}^T)^j \nabla u^{3-i}_s+v^{3-i,j}_s,u^i_s)\ge-%
\frac{1}{4}|\|(k^+k\bar{g}^T)^j\nabla u^{3-i}_s+v^{3-i,j}_s|\|^2 -C\|
u^i_s\|^2\,.
\end{equation*}
Moreover 
\begin{align*}
(\alpha^{j}v^{j},u^i)=(v^{i,j}_s+(k^+k\bar{g}^T)^j\nabla u^i_s,
\alpha^{j}u^i_s)-((k^+k\bar{g}^T)^j\nabla u^i_s, \alpha^{j}u^i_s).
\end{align*}
We treat the first term on the right as in the last inequality, and by
integration by parts the second term is bounded from below by $-C\|u^i_s\|^2$%
. Now assuming that we can take the expectation in \eqref{energyequ} and
that the expectation of the stochastic integral vanishes (this is not a
serious difficulty, although it requires the use of stopping times), and
combining the resulting identity with \eqref{ineq} and Gronwall's Lemma, we
deduce \eqref{estim1} with $m=0$, at least with the $\sup_t$ outside the
expectation. \eqref{estim1} then follows using Doob's inequality. The reader
may have noticed that the above argument requires that all coefficients have
bounded first order partial derivatives, and the entries of the matrix $a$
have also bounded second order partial derivatives.

Next we approximate our pair of BSPDEs with a system indexed by $%
\varepsilon>0 $, where the operator $A$ has been replaced by $%
A_\varepsilon:=A+\varepsilon\Delta$, where $\Delta$ stands for the Laplace
operator. We can now invoke Theorem \ref{th:SP} to obtain the existence of a
solution $(u^1_\varepsilon,u^2_\varepsilon)$ to our approximate system of
BSPDEs. Clearly the above computations yield that $(u^i_\varepsilon,v^i_%
\varepsilon)$, $i=1,2$ satisfies the estimate \eqref{estim1} with $m=0$ and
a constant $C$ which is independent of $\varepsilon$. Hence we can extract a
subsequence such that each pair $(u^i_\varepsilon,v^i_\varepsilon)$
converges weakly $L^2_{\mathcal{P}}((0,T)\times\Omega;H^0\times(H^0)^{%
\otimes\ell^{\prime }})$, and it is not too hard to show that the limit
still satisfies \eqref{estim1} for $m=0$, and solves the system of BSPDEs
(we take the limit in the equation written in weak form).

However, we are interested in more regular solutions. Mimicking the computations
done in \cite{DTZ}, we can extend the above computations to estimate the
norms in $H^m$. Of course, this must be done sequentially w.r.t. $m$. It
requires to take partial derivatives in our system of BSPDEs, which
introduces a forcing term involving lower order derivatives, but it can be
handled in our situation similarly as both in \cite{DTZ} and in \cite{DM}.
As a result, we can take the limit weakly in $L^2_{\mathcal{P}%
}((0,T)\times\Omega;H^m\times(H^m)^{\otimes\ell^{\prime }})$. The weak
continuity of $u^i$ with values in $H^m$ follows by standard arguments. The
result follows.
\end{proof}

\begin{remark}
Theorem 2.1 in \cite{DTZ} also asserts the uniqueness of the corresponding
equation. Whilst we don't need it for our purpose, the uniqueness of the
solution of (\ref{splusf}) can be shown in a similar manner and it is a
consequence of \eqref{estim1}.
\end{remark}

\begin{theorem}
\label{th:systemofBSPDEs} Let $\theta _{t}$ be the $\mathbb{C}$--valued
solution of the SDE 
\begin{equation}  \label{phi}
d\theta _{t}=i\theta _{t}r_{t}^{j}d\tilde{W}_{t}^{j},\quad \theta _{0}=1\,,
\end{equation}%
where $r_{t}$ is an arbitrary element of $L^{\infty }([0,T];\mathbb{R}%
^{d^{\prime }})$. Under the 
assumptions \textbf{E} and \textbf{U}, we get that, $(u,v)$ denoting the
solution of the BSPDE \eqref{SBSPDE} and $\pi_t$ a solution of the Zakai
equation \eqref{eq:filterEqns:zakaifunctint} which satisfies 
\begin{equation}  \label{hypsolZak}
\mathbb{E}\left[\sup_{0\le t\le T}\pi_t(\mathbf{1})^2\right],
\end{equation}
we have 
\begin{equation}  \label{ii}
d\theta _{t}\pi _{t}(u _{t})=\theta _{t}\pi _{t}(B^{j}u
_{t}+v_{t}^{j}+ir_{t}^{j}u _{t})d\tilde{W}_{t}^{j}.
\end{equation}%
Moreover the process $\{\theta _{t}\pi _{t}(u _{t}),\ 0\le t\le T\}$ is a
martingale, and therefore $\tilde{\mathbb{E}}\left[ \theta _{T}\pi
_{T}^{1}(\varphi )\right] =\tilde{\mathbb{E}}\left[ \pi _{0}^{1}(\psi _{0})%
\right] $.
\end{theorem}

\begin{proof}
Thanks to our assumption, it follows from Theorem \ref{bspde} that we can
apply Theorem \ref{Itoform}, from which we deduce that 
\begin{align}
d\pi _{t}(u_{t})& =\pi _{t}\left( Au _{t}-Au
_{t}-B^{j}v_{t}^{j}-ir_{t}^{j}B^{j}u
_{t}-ir_{t}^{j}v_{t}^{j}+B^{j}v_{t}^{j}\right) dt+\pi _{t}(B^{j}u
_{t}+v_{t}^{j})d\tilde{W}_{t}^{j}  \notag \\
& =-ir_{t}^{j}\pi _{t}\left( B^{j}u_{t}+v_{t}^{j}\right) dt+\pi
_{t}(B^{j}u_{t}+v_{t}^{j})d\tilde{W}_{t}^{j}  \label{i}
\end{align}%
The identity \eqref{ii} follows from \eqref{phi} and \eqref{i} by It\^o's
chain rule. The martingale property of the process $\{\theta _{t}\pi _{t}(u
_{t}),\ 0\le t\le T\}$ will follow from the Burkholder--Davis--Gundy
inequality for the first moment and the following estimate 
\begin{equation}  \label{Ubound}
\tilde{\mathbb{E}}\left[\sqrt{\int_0^T |\theta _{t}\pi
_{t}(B^{j}u_{t}+v_{t}^{j}+ir_{t}^{j}u _{t})|^2ds} \right]<\infty
\end{equation}
for $j=1,...,\ell^{\prime }$ which we now establish. By Theorem \ref{bspde},
our current assumptions imply that for some $m>d/2+2$, $u\in L^2(\Omega\times%
[0,T];H^{m}(\mathbb{R}^d))$ and for $1\le j\le \ell^{\prime }$, $v^j\in
L^2(\Omega\times[0,T];H^{m-1}(\mathbb{R}^d))$. Consequently, thanks to the
Sobolev embedding theorem, we deduce that 
\begin{equation*}
\tilde{\mathbb{E}}\int_0^T\left[\sup_x|\nabla
u(t,x)|^2+\sup_x|u(t,x)|^2+\sum_{j=1}^m\sup_x|v^j(t,x)|^2\right]dt<\infty\,.
\end{equation*}
So, if we define $C^j(t,.):=B^ju_t+v^j_t+ir^j_tu_t$, we have that for $1\le
j\le \ell^{\prime }$, 
\begin{equation}  \label{smoothness}
\tilde{\mathbb{E}}\int_0^T\sup_{x}|C^j(t,x)|^2dt<\infty\,.
\end{equation}
We first note that for any $0\le t\le T$, 
\begin{equation*}
|\theta_t|=\exp\left(\frac{1}{2}\sum_{j=1}^{\ell^{\prime }}\int_0^t|r^j_s|^2
ds\right)
\end{equation*}
is a deterministic quantity, whose supremum over $0\le t\le T$ is finite.
Hence we have 
\begin{align*}
\int_0^T|\theta_t\pi_t(C^j(t,\cdot))|^2dt&\le \sup_{0\le t\le
T}|\theta_t|^2\sup_{0\le t\le T}\pi_t(\mathbf{1})^2\int_0^T%
\sup_x|C^j(t,x)|^2dt, \\
\tilde{\mathbb{E}}\left(\sqrt{\int_0^T|\theta_t\pi(C^j(t,\cdot))|^2dt}%
\right) &\le\sup_{0\le t\le T}|\theta_t|\tilde{\mathbb{E}}\left\{\sup_{0\le
t\le T}\pi_t(\mathbf{1}) \sqrt{\int_0^T\sup_x|C^j(t,x)|^2dt}\right\} \\
&\le \sup_{0\le t\le T}|\theta_t|\ \sqrt{\tilde{\mathbb{E}}\left(\sup_{0\le
t\le T}\pi_t(\mathbf{1})^2\right)} \sqrt{\tilde{\mathbb{E}}%
\int_0^T\sup_x|C^j(t,x)|^2dt},
\end{align*}
and \eqref{Ubound} now follows from \eqref{hypsolZak} and \eqref{smoothness}.
\end{proof}

%
%

\begin{theorem}
\label{th:Unique} Under assumptions {\bf E} and {\bf U}, there exists a unique solution
of the equation \eqref
{eq:filterEqns:zakaifunctint} in the class of $\mathcal{P}(\mathcal{Y}_t)$-measurable
measure valued processes satisfying 
\eqref{hypsolZak} for any $T>0$.
Since the uniqueness of the solution of equation (\ref{eq:filterEqns:zakaifunctint}) is equivalent to that of 
equation (\ref{eq:filterEqns:ks}) following from Theorem \ref{equivuniq}, we also deduce the uniqueness of the solution of equation (\ref{eq:filterEqns:ks}).

\end{theorem}

\begin{proof}
Assume that there are two solutions of the equation (\ref%
{eq:filterEqns:zakaifunctint}) denoted by $\pi _{1},\pi _{2}$. We observe
the following sequence of identities 
\begin{equation*}
\tilde{\mathbb{E}}\left[ \theta _{T}\pi _{T}^{1}(\varphi )\right] =\tilde{ 
\mathbb{E}}\left[ \theta _{0}\pi _{0}^{1}(u_{0})\right] =\tilde{\mathbb{E}}%
\left[ \theta _{0}\pi _{0}^{2}(u_{0})\right] =\tilde{\mathbb{E}}\left[
\theta _{T}\pi _{T}^{2}(\varphi )\right]
\end{equation*}%
and since both $S_{T}$ is a total set, and $\varphi$ is an arbitrary smooth
function, it follows that $\pi _{T}^{1}=\pi _{T}^{2}$.
\end{proof}

Finally we note that the unnormalized conditional distribution satisfies the
condition \eqref{hypsolZak}. Indeed, if we let now $\pi _{t}$ denote that
unnormalized conditional distribution at time $t$, we deduce from Lemma \ref%
{le:filterEqns:rhoCadlag} that $\pi _{t}(\mathbf{1})=\tilde{Z}_{t}$, which
is a $\tilde{\P }$ martingale, hence from Doob's inequality, it suffices to
prove that $\tilde{\mathbb{E}}[|\tilde{Z}_{t}|^{2}]<\infty $ for all $t>0$,
which follows from \eqref{ztilde2} and the boundedness of $h_{2}$.

\section{Appendix}
In this Appendix, we recall the definition of the
 Moore--Penrose pseudo inverse of a possibly rectangular matrix, and prove that the map which to a matrix associates its pseudo--inverse is measurable.

In what follows, $A^{+}\in \mathbb{R}^{\ell\times d}$
stands for the Moore-Penrose pseudo-inverse of the matrix $A\in\R^{d\times\ell}$. The
Moore-Penrose pseudo-inverse $A^{+}$ of the matrix $A$ is uniquely characterised
by the following four properties 
\begin{equation}
AA^{+}A=A,~~~A^{+}AA^{+}=A^{+},~~\left( A^{+}A\right) ^{\top
}=A^{+}A,~~~~~\left( AA^{+}\right) ^{\top }=AA^{+},  \label{moorepenrose}
\end{equation}
see \cite{bengrev} for details.

Moreover from the identities in (\ref{moorepenrose}) we deduce that $A^{+}A$ is
a projection onto the range of $A^{+}A$ and $I-A^{+}A$ is a projection onto
the orthogonal space of the range of $A^{+}A$. In particular we deduce that the norm of
 $A^{+}A$ as a linear operator is bounded by $1$. Since all norms on a finite dimensional space are equivalent, for any matrix norm $\|\cdot\|$, there exists a constant $C$ which depends only upon $\ell$ and the particular choice of a norm on the set of $\ell\times \ell$ matrices such that 
 \begin{equation}\label{k+k}
 \|A^+A\|\le C,
 \end{equation} for any positive integers $d,l$ and any $d\times\ell$ matrix $A$.

Next we give details of an explicit construction of the pseudo inverse of a
matrix that will help us prove the measurability of the mapping $A\mapsto A^+$. In what
follows we use the notation $A$ for a generic matrix $A\in \mathbb{R}^{d\times \ell }$. We follow here the construction and the analysis in \cite{bengrev}.

Let $0\leq r\leq \min \left( d,\ell \right) $ be the
rank of $A$. If $r=0$ this means that $A=O^{d\times \ell}$
where $O^{d\times \ell}$ is the matrix with all entries
null. In this degenerate case, $A^{+}=O^{\ell\times d}$ is
the matrix with all null entries. If $r>0$, then $A$ has a invertible minor
of order $r$. Following Theorem 5 page 48 in \cite{bengrev}, the
pseudo-inverse is given by 
\begin{equation*}
A^{+}=G^{\top }\left( F^{\top }AG^{\top }\right) ^{-1}F^{\top },
\end{equation*}
where $G$ and $F$ are matrices which appear in a full-rank
factorization of $A$: 
\begin{equation*}
A=FG,~~F\in \mathbb{R}^{d\times r},~~~G\in \mathbb{R}^{r\times\ell }.
\end{equation*}
The pair $(F,G)$ is not unique, but we shall give one construction which is based upon a particular choice of an invertible $r\times r$ minor of $A$.
We describe briefly the construction of a full-rank factorization of $A$
 (see also Section 4, in particular page 26 in \cite{bengrev}): 
 
 Let $1\le c_1 < c_2 < \cdots < c_r \le d\wedge\ell$ denote the ranks of the $r$ columns containing all terms of one arbitrarily chosen
 invertible $r\times r$ minors of $A$. Let now $P$ be the permutation matrix, which when applied to $A$ on the right, makes the $c_i$--th
 column of $A$ into the $i$--th column of $AP$, $1\le i\le r$. Let next $P_1$ denote the submatrix of $P$ consisting of its first $r$ columns, and 
 $F=AP_1$. Finally let $G$ be the unique $r\times \ell$ matrix such that $A=F G$. It is clear (see Lemma 1 page 26 of \cite{bengrev}) that the terms of the $i$--th column of $G$ 
are the unique coefficients which express the $i$--th column of A as a linear combination of the elements of the basis of $R(A)$ given by the columns of $F$.
We have $A=FG$, where $F$ (resp. $G$) is a $d\times r$ (resp. $r\times \ell$) matrix of rank $r$.

We now want to prove the following lemma.

\begin{lemma}\label{MP-measurable} 
The mapping $A\mapsto A^+$ is measurable from $\R^{d\times\ell}$ into $\R^{\ell\times d}$.
\end{lemma}
We note that the above mapping is clearly not continuous (in the case $d=\ell=1$, for $A\in\R$, $A^+=1/A$ if $A\not=0$, and $0^+=0$,
so if $A_n>0$, $A_n\to0$, then $A_n^+\to+\infty$, while $(\lim_n A_n)^+=0$).

\begin{proof}
We first need to find a consistent way of
identifying an invertible minor of order $r$ of the matrix $A$. For
this we introduce the following enumeration of the minors. 
The matrix $A$ has $d\ell $ minors or order 1, $\left( 
\begin{array}{c}
d \\ 
2
\end{array}
\right) \left( 
\begin{array}{c}
\ell  \\ 
2
\end{array}
\right) $ minors of order 2, ..., and $\left( 
\begin{array}{c}
\max(d, \ell) \\ 
\min(d, \ell)
\end{array}
\right) $ minors of order $\min(d, \ell)$. By convention we add a `minor' of order 0 (to account for the
case $r=0$) whose determinant is chosen to be $0$. \ Next we consider 
\begin{equation*}
S\left(A\right) \in \mathbb{R}^{\theta },~~~~\theta :=1+d\ell
+\left( 
\begin{array}{c}
d \\ 
2
\end{array}
\right) \left( 
\begin{array}{c}
\ell  \\ 
2
\end{array}
\right) +...+\left( 
\begin{array}{c}
\max(d, \ell) \\ 
\min(d, \ell)
\end{array}
\right) ,
\end{equation*}
\begin{equation*}
S\left( A\right) : =\left( 0,....\right) ,
\end{equation*}
where $S(A)$ is the list of the determinants of all the minors of $A$, and we choose a fixed arbitrary enumeration of  all those minors.
If $r>0$, we denote by $m\left( A\right) $ the highest index in the
enumeration of the minors for which the corresponding determinant is
non-zero. Such an index exists and $m\left( A\right) >1$. This means that 
\begin{equation*}
S\left( A\right) =\left( 0,...,d\left( A\right) ,0,...,0\right) ,
\end{equation*}
where $d(A)$ is the $m\left( A\right)$-th entry of the vector $S\left(
A\right) ,$ that is $d(A)=S\left( A\right) _{m\left( A\right) }$ is
not zero and all subsequent entries (if any) $S\left( A\right) _{m\left(
A\right) +1}$, $S\left( A\right) _{m\left( A\right) +2},...$ are zero.

If $r=0$,  then $S\left( A\right)$ has all entries null 
\begin{equation*}
S\left( A\right) =\left( 0,...,...,0\right), 
\end{equation*}
and $m\left( A\right) =1$, by convention.

 Note that the function $A\mapsto
m\left( A\right) $ is an integer valued function.
We split $\mathbb{R}^{\theta }$ into a finite collection of disjoint sets $
\left\{ H_{n}\right\} _{n=1}^{\theta }$ such that the index $m\left(
A\right) $ stays constant if $S\left(A\right) $ takes values in $H_{n}$ :

\begin{itemize}
\item $H_{1}\subset \mathbb{R}^{\theta },~~H_{1}=\left\{
(0,0,0,...,0)\in \mathbb{R}^{\theta }\right\} $. On this set $m\left(
A\right) =1$. In this case, since $S\left(A\right) =(0,0,...,0)$, $A$ is the null matrix,

\item $H_{2}\subset \mathbb{R}^{~\theta },~~H_{2}=\left\{
(0,a_{2},0,...,0)\in \mathbb{R}^{~\theta },~a_{2}\neq 0\right\} $. On
this set $m\left( A\right) =2$. In this case, the first ranked minor of $A$ 
is not zero and the determinants of all the higher ranked minors are zero.

\item $H_{3}\subset \mathbb{R}^{~\theta },~~H_{3}=\left\{
(0,a_{2},a_{3},...,0)\in \mathbb{R}^{~\theta },~a_{3}\neq 0\right\} $.
On this set $m\left(A\right) =3$. In this case, the second ranked minor of the matrix $A$ 
 is not zero and the determinants of all the higher ranked minors are zero.

\item ....

\item $H_{\theta }\subset \mathbb{R}^{~\theta },~~H_{\theta }=\left\{
(0,a_{2},a_{3},...,a_{\theta })\in \mathbb{R}^{~\theta },~a_{\theta
}\neq 0\right\} $. On this set $m\left(A\right) =\theta $. In this case,
the last minor of $A$ in the list has a non zero determinant.
\end{itemize}

We distinguish two cases:

If $m\left( A\right) $ takes values in the set $H_{1}$ (in other words $
m\left( A\right) =1$ and $S\left( A\right) =
(0,0,...,0)$ ), then $A$ is the null matrix $O^{d\times \ell }$ (in other words it is constant) on the preimage of
this set%
\begin{equation*}
Q_{1}=\left\{ A \in \mathbb{R}^{d\times\ell}|S\left(A\right) =(0,0,...,0)\in H_{1}\right\} .
\end{equation*}

If $m\left( A\right) $ takes values in each of the remaining sets $%
H_{i},i=2,...,\theta $, then $A$ has one fixed invertible
minor on the preimage of each of those sets. 
\begin{equation*}
Q_{i}=\left\{ A\in \mathbb{R}^{d\times\ell}|S\left( A\right) \in H_{i}\right\} ,~~~i=2,3,...,\theta .
\end{equation*}%
As a result, on each set $Q_{i}$, the same (in most cases arbitrarily chosen) invertible $r\times r$ minor is selected. Let $P_i$ be the permutation matrix presented above which, when 
applied to $A$ on the right, moves the $r$ columns containing elements of the selected invertible minor into the first $r$ columns. It is clear that those columns constitute a basis of 
$R(A)$, the range of $A$. Denote by $P_{i,1}$ the $\ell\times r$ matrix consisting of the first $r$ columns of $P_i$. The matrix $P_{i,1}$ is constant on $Q_{i}$, hence $F_i=AP_{i,1}$ is a continuous function of $A$ on $Q_{i}$. 
Next let $G_i$ be the unique $r\times \ell$ matrix such that $A=F_i G_i$.  
The
matrix $G_i$ is a rational, hence continuous function of the entries of $A$.
Therefore on each set $Q_{i}$ \ the Moore-Penrose pseudo-inverse of $A$ 
 \begin{equation*}
A^+=G_i^{\top } \left( F_i^{\top }A G_i^{\top } \right) ^{-1}F_i^{\top }
\end{equation*}
is a continuous function of $A$.
Finally, since $\left\{ A \in Q_{i}\right\} $ is Borel, we
get that 
\begin{eqnarray*}
A^+ &=&O^{d\times \ell }1_{A \in Q_{1} } +\sum_{i=2}^{\theta }G_i^{\top } \left( F_i^{\top } A G_i^{\top } \right)^{-1}F_i^{\top } 1_{\left\{ A \in Q_{i}\right\} }
\end{eqnarray*}
is a measurable function of $A$.
\end{proof}

\end{document}